\documentclass[11pt,leqno]{article}

\usepackage[english]{babel}

\usepackage{amsmath}

\usepackage{amssymb}

\usepackage{xypic}

\usepackage[all]{xy}

\usepackage{latexsym}

\usepackage{amsfonts}

\usepackage{epsf}

\usepackage[dvips]{graphicx}

\usepackage[latin1]{inputenc}

\usepackage[T1]{fontenc}

\numberwithin{equation}{section}

\newcommand{\C}{\mathcal{C}}

\newcommand{\D}{\mathcal{D}}

\newcommand{\M}{\mathcal{M}}

\newcommand{\W}{\mathcal{W}}

\renewcommand{\mod}{\mathrm{Mod}}

\newcommand{\ind}{\mathrm{Ind}}

\newcommand{\cov}{\mathrm{Cov}}

\newcommand{\coh}{\mathrm{Coh}}

\newcommand{\op}{\mathrm{Op}}

\newcommand{\rc}{\mathbb{R}\textrm{-}\mathrm{c}}

\newcommand{\CC}{\mathbb{C}}

\newcommand{\R}{\mathbb{R}}

\newcommand{\Z}{\mathbb{Z}}

\newcommand{\N}{\mathbb{N}}

\newcommand{\OO}{\mathcal{O}}

\newcommand{\NN}{\mathcal{N}}

\newcommand{\LL}{\mathcal{L}}

\newcommand{\II}{\mathcal{I}}

\newcommand{\ltens}{\overset{L}{\otimes}}

\newcommand{\osi}{\stackrel{\sim}{\gets}}

\newcommand{\iso}{\stackrel{\sim}{\to}}

\newcommand{\dbtxr}{\mathcal{D}\mathit{b}^t_{X_\mathbb{R}}}

\newcommand{\dbt}{\mathcal{D}\mathit{b}^t}

\newcommand{\db}{\mathcal{D}\mathit{b}}

\newcommand{\ot}{\mathcal{O}^t}

\newcommand{\ddxy}{\mathcal{D}_{\xmenoy}}

\newcommand{\ddyx}{\mathcal{D}_{\ymenox}}

\newcommand{\xmenoy}{{X}\rightarrow{Y}}

\newcommand{\ymenox}{{Y}\leftarrow{X}}

\newcommand{\tho}{\mathit{T}\mathcal{H}\mathit{om}}

\newcommand{\Th}{\mathrm{THom}}

\newcommand{\OW}{\OO^\mathrm{w}}
\newcommand{\OWX}{\OO^\mathrm{w}_X}
\newcommand{\OWY}{\OO^\mathrm{w}_Y}

\newcommand{\CWM}{\C^{{\infty ,\mathrm{w}}}_M}

\newcommand{\CW}{\C^{{\infty ,\mathrm{w}}}}
\newcommand{\wtens}{\overset{\mathrm{w}}{\otimes}}

\newcommand{\sol}{\mathcal{S}ol}

\newcommand{\rh}{\mathit{R}\mathcal{H}\mathit{om}}

\newcommand{\ho}{\mathcal{H}\mathit{om}}

\newcommand{\Ho}{\mathrm{Hom}}

\newcommand{\Rh}{\mathrm{RHom}}

\newcommand{\id}{\mathrm{id}}

\renewcommand{\dim}{\textbf{Proof.}}

\newcommand{\qed}{\nopagebreak \phantom{} \hfill $\Box$ \\}

\newcommand{\pE}{\overset{\hspace{0.1cm}}{\dot{E}}}

\newcommand{\pU}{\overset{\hspace{0.1cm}}{\dot{U}}}

\newcommand{\pV}{\overset{\hspace{0.1cm}}{\dot{V}}}

\newcommand{\supp}{\mathrm{supp}}

\newcommand{\pittau}{\stackrel{\mathbf{.}}{\tau}}

\newcommand{\RP}{\mathbb{R}^{{\scriptscriptstyle{+}}}}

\newcommand{\imin}[1]{#1^{-1}}

\newcommand{\lind}[1]{\underset{#1}{\underrightarrow{\lim}}}

\newcommand{\Lind}{\underrightarrow{\lim}}  

\newcommand{\lpro}[1]{\underset{#1}{\underleftarrow{\lim}}}

\newcommand{\exs}[3]{0 \to {#1} \to {#2} \to {#3} \to 0}

\newcommand{\lexs}[3]{0 \to {#1} \to {#2} \to {#3}}

\newcommand{\dt}[3]{{#1} \to {#2} \to {#3} \stackrel{+}{\to}}

\newtheorem{teo}{Theorem}[subsection]

\newtheorem{df}[teo]{Definition}

\newtheorem{cor}[teo]{Corollary}

\newtheorem{oss}[teo]{Remark}

\newtheorem{prop}[teo]{Proposition}

\newtheorem{lem}[teo]{Lemma}

\newtheorem{es}[teo]{Exemple}

\newtheorem{nt}[teo]{Notations}

\author{Luca Prelli}
\title{\bf{CONIC SHEAVES ON SUBANALYTIC SITES AND LAPLACE TRANSFORM}}
\date{}

\usepackage{latexsym}
\begin{document}


\maketitle




\thispagestyle{empty}

\begin{abstract} Let $E$ be a $n$ dimensional complex vector space and
let $E^*$ be its dual. We construct the conic sheaves
$\OO^t_{E_{\RP}}$ and $\OO^{\mathrm{w}}_{E_{\RP}}$ of tempered and
Whitney holomorphic functions respectively and we give a sheaf theoretical interpretation of the
Laplace isomorphisms of \cite{KS97} which gives the isomorphisms in the
derived category $\OO^{t\land}_{E_{\RP}}[n] \simeq
\OO^t_{E^*_{\RP}}$ and $\OO^{\mathrm{w}\land}_{E_{\RP}}[n] \simeq
\OO^{\mathrm{w}}_{E^*_{\RP}}$.
\end{abstract}


\tableofcontents

\addcontentsline{toc}{section}{\textbf{Introduction}}

\section*{Introduction}

Classical sheaf theory is not well suited to the study of various
objects in Analysis, which are not defined by local properties. To
overcome this problem Kashiwara and Schapira in \cite{KS01}
developed the theory of ind-sheaves on a locally compact space and
defined the six Grothendieck operations in this framework. For a
real analytic manifold $X$, they defined the subanalytic site
$X_{sa}$ as the site generated by the category of subanalytic open
subsets and whose coverings are the locally finite coverings in
$X$. Moreover they proved the equivalence between the categories
of ind-$\R$-constructible sheaves and sheaves on $X_{sa}$. They
constructed the sheaf of tempered distributions as a sheaf on
$X_{sa}$, and when $X$ is a complex manifold they introduced the
sheaf of tempered holomorphic functions.
Thanks to the results of \cite{Pr1} we have a direct construction
of the six Grothendieck operations on $\mod(k_{X_{sa}})$, without
using the theory of ind-sheaves. So we will work directly with
sheaves on subanalytic sites.

Let $X$ be a real analytic manifold
endowed with an action of $\RP$. Our first goal is to define conic sheaves on $X_{sa}$ and then, when
$X$ is a vector bundle, to extend the construction of the
Fourier-Sato transform for classical sheaves to conic sheaves on
subanalytic sites. In order to do that we have to choose a suitable definition of conic sheaf: indeed there are several definitions, which are equivalent in the classical case but not in the framework of subanalytic sheaves. We choose the one which satisfies some desirable properties (as the equivalence with sheaves on the conic topology associated to the action and the isomorphism of conic sheaves with limits of conic $\R$-constructible sheaves) and for which the Fourier-Sato isomorphism applies.

Let $E$ be a complex vector space. As an application we construct the conic
sheaves $\OO^t_{E_{\RP}}$ and $\OO^{\mathrm{w}}_{E_{\RP}}$ of tempered
and Whitney holomorphic functions respectively and we
prove a sheaf theoretical interpretation of the Laplace isomorphisms of \cite{KS97} which induce
isomorphisms in the derived category
$\OO^{t\land}_{E_{\RP}}[n] \simeq \OO^t_{E^*_{\RP}}$ and
$\OO^{\mathrm{w}\land}_{E_{\RP}}[n] \simeq \OO^{\mathrm{w}}_{E^*_{\RP}}$.\\

In more details, the contents of this paper are as follows.

In Section 1 we first recall the results of \cite{KS01} and
\cite{Pr1} on sheaves on su\-ba\-na\-ly\-tic sites. Then we consider the category of conic sheaves on subanalytic sites.
Let $k$ be a field and let $X$ be a real analytic manifold with an action $\mu$ of $\RP$. Let $U$ be an open subset of $X$. We say that $U$
is $\RP$-connected if its intersections with the orbits of $\mu$
are connected. We denote $\RP U$ the conic open set associated to
$U$ (i.e. $\RP U=\mu(U,\RP)).$
A sheaf $F$ on $X_{sa}$ is
said to be conic if $\Gamma(\RP U;F) \iso \Gamma(U;F)$ for each
$\RP$-connected relatively compact open subanalytic subset $U$ of
$X$. We call $\mod_{\RP}(k_{X_{sa}}) \subset \mod(k_{X_{sa}})$
the category of conic sheaves. This definition is different from the classical one.
Let us consider the projection $p:X \times \RP \to X$. One can define the subcategory
$\mod^\mu(k_{X_{sa}})$ of $\mod(k_{X_{sa}})$ consisting of sheaves
satisfying $\imin \mu F \simeq \imin p F$. The categories
$\mod^\mu(k_{X_{sa}})$ and $\mod_{\RP}(k_{X_{sa}})$ are not
equivalent in general.
The category of conic subanalytic sheaves has many good
properties, for example it is equivalent to the category of
sheaves on the conic subanalytic site $X_{sa,\RP}$ (i.e. the
category of open conic subanalytic subsets of $X$ with the
topology induced from $X_{sa})$. This equivalence is strictly
related to the geometry of subanalytic open subsets of $X$.
When $E$ is a vector bundle, one can define the
Fourier-Sato transform which gives an equivalence between conic
sheaves on $E_{sa}$ and conic sheaves on $E^*_{sa}$, where $E^*$
denotes the dual vector bundle.

In Section 2 we study the conic sheaves of tempered and Whitney holomorphic
functions. Let $E$ be a real vector space, and let $E
\hookrightarrow P$ be its projective compactification. We define the
conic sheaves of tempered distributions $\dbt_{E_{\RP}}$ and Whitney
$\C^\infty$-functions $\CW_{E_{\RP}}$. If $U$ is an open subanalytic
cone, the sections $\Gamma(U;\dbt_{E_{\RP}})$ are distributions which
are tempered on the boundary of $U$ and at infinity, and the sections
$\Gamma(U;\CW_{E_{\RP}})$  are Whitney $\C^\infty$-functions on $U$
with rapidly decay at infinity.
If $F$ is a conic
$\R$-constructible sheaf on $E$ we have
\begin{eqnarray*}
\Rh(F,\dbt_{E_{\RP}}) & \simeq & \Th(F,\db_E), \\
\Rh(F,\CW_{E_{\RP}}) & \simeq
& D'F \stackrel{\rm W}{\otimes}\C^\infty_E,
\end{eqnarray*}
where $\Th(\cdot,\db_E)$ and $\cdot
\stackrel{\rm W}{\otimes}\C^\infty_E$ are the functors introduced in \cite{KS97}.

 When $E$ is an $n$ dimensional complex vector space we
define the sheaves $\ot_{E_{\RP}}$ and $\OW_{E_{\RP}}$ of tempered and Whitney holomorphic functions
taking the solutions of the Cauchy-Riemann system with values in
tempered distributions and Whitney $\C^\infty$-functions respectively.
 We show that these sheaves are invariant by the Laplace
transform. In fact the Laplace isomorphisms of \cite{KS97} induce
isomorphisms  in the derived category $\OO^{t\land}_{E_{\RP}}[n]
\simeq \OO^t_{E^*_{\RP}}$ and $\OO^{\mathrm{w}\land}_{E_{\RP}}[n]
\simeq \OO^{\mathrm{w}}_{E^*_{\RP}}$, where $\land$ denotes the
extension of the Fourier-Sato transform to conic sheaves on
$E_{sa}$. Moreover these isomorphisms are compatible with the action of the Weyl
algebra.

\section{Conic sheaves on subanalytic sites}

In the following $X$ will be a real analytic manifold and $k$ a
field. Reference are made to
\cite{KS01,Pr1} for an introduction to sheaves on subanalytic sites. We refer to \cite{BM88,Lo93} for the theory of subanalytic
sets and to \cite{Co00,VD98} for the more general theory of o-minimal structures.

\subsection{Review on sheaves on subanalytic sites}


Denote by $\op(X_{sa})$ the category of subanalytic open subsets of
$X$. One endows $\op(X_{sa})$ with the following topology: $S
\subset \op(X_{sa})$ is a covering of $U \in \op(X_{sa})$ if for
any compact $K$ of $X$ there exists a finite subset $S_0\subset S$
such that $K \cap \bigcup_{V \in S_0}V=K \cap U$. We will call
$X_{sa}$ the subanalytic site.\\

Let $\mod(k_{X_{sa}})$ denote the category of sheaves on $X_{sa}$.
Then $\mod(k_{X_{sa}})$ is a Grothendieck category, i.e. it admits
a generator and small inductive limits, and small filtrant
inductive limits are exact. In particular as a Grothendieck
category, $\mod(k_{X_{sa}})$ has enough injective objects.\\

Let $\mod_{\rc}(k_X)$  be the abelian category of
$\R$-constructible sheaves on $X$, and consider its subcategory
$\mod^{c}_{\rc}(k_X)$ consisting of sheaves whose support is
compact.\\

We denote by $\rho: X \to X_{sa}$ the natural morphism of sites.
We have functors


\begin{equation*}
\xymatrix{\mod(k_X)   \ar@ <5pt> [rr]^{\rho_*} \ar@
<-5pt>[rr]_{\rho_!} &&
  \mod(k_{X_{sa}}) \ar@ <0pt> [ll]|{\imin \rho}. }
\end{equation*}

The functors $\imin \rho$ and $\rho_*$ are the functors of inverse
image and direct image respectively. The sheaf $\rho_!F$ is the
sheaf associated to the presheaf $\op_{sa}(X) \ni U \mapsto
F(\overline{U})$. In particular, for $U \in \op(X)$ one has
$\rho_!k_U \simeq \lind {V \subset \subset U} \rho_*k_V$, where $V
\in \op_{sa}(X)$. Let us summarize the properties of these
functors:

\begin{itemize}

\item the functor $\rho_*$ is fully faithful and left exact, the restriction
of $\rho_*$ to $\mod_{\rc}(k_X)$ is exact,
\item the functor $\imin \rho$ is exact,
\item the functor $\rho_!$ is fully faithful and exact,
\item $(\imin \rho,\rho_*)$ and $(\rho_!,\imin \rho)$ are pairs of
adjoint functors.
\end{itemize}

For each $F \in \mod(k_{X_{sa}})$, there exists a small filtrant
inductive system $\{F_i\}$, with $F_i \in \mod^{c}_{\rc}(k_X)$,
such that $F \simeq \lind i \rho_*F_i$. Moreover let $\{F_i\}$ be
a filtrant inductive system of $k_{X_{sa}}$-modules and let $G \in
\mod^c_{\rc}(k_X)$. One has the isomorphism
\begin{equation}
\Ho(\rho_*G,\lind iF_i) \simeq \lind i \Ho(\rho_*G,F_i).
\end{equation}

Let $X,Y$ be two real analytic manifolds, and let $f:X \to Y$ be a
real analytic map. We have a commutative diagram
\begin{equation}
\xymatrix{X \ar[d]^\rho \ar[r]^f & Y \ar[d]^\rho \\
X_{sa} \ar[r]^f & Y_{sa}}
\end{equation}

We get external operations $\imin f$, $f_*$ and $f_{!!}$, where
the notation $f_{!!}$ follows from the fact that $f_{!!} \circ
\rho_* \not\simeq \rho_* \circ f_!$ in general. If $f$ is proper on
$\supp(F)$ then $f_*F \simeq f_{!!}F$, in this case $f_{!!}$
commutes with $\rho_*$. While functors $\imin f$ and $\otimes$ are
exact, the functors $\ho$, $f_*$ and $f_{!!}$ are left exact and
admit right derived functors. In particular the functor $Rf_{!!}$
admits a right adjoint, denoted by $f^!$, and we get the usual
isomorphisms between Grothendieck operations (projection formula,
base change formula, K\"unneth formula, etc.) in the framework of
subanalytic sites.
We refer to \cite{Pr1} for a
detailed exposition.

Finally we recall the relations between the six Grothendieck
operations and the functors $\imin \rho$, $R\rho_*$ and $\rho_!$.

\begin{itemize}
\item the functor $\imin \rho$ commutes with $\otimes$, $\imin f$ and
$Rf_{!!}$,
\item the functor $R\rho_*$ commutes with $\rh$, $Rf_*$ and
$f^!$,
\item the functor $\rho_!$ commutes with $\otimes$ and $\imin f$,

\item if $f$ is a topological submersion (i.e. it is locally
isomorphic to a projection $Y \times \R^n \to Y$), then $f^!\simeq
\imin f \otimes f^!k_Y$ commutes with $\imin \rho$ and $Rf_{!!}$
commutes with $\rho_!$.
\end{itemize}

\subsection{Conic subanalytic sheaves}

Let $X$ be a real analytic manifold endowed with an analytic
action $\mu$ of $\RP$. We have a diagram

$$\xymatrix{X \ar[r]^{\hspace{-5mm}\iota} & X \times \RP \ar@ <2pt>
[r]^{\hspace{0.3cm}\mu} \ar@ <-2pt> [r]_{\hspace{0.3cm}p}& X,}$$
where $\iota(x)=(x,1)$ and $p$ denotes the projection. We have $\mu
\circ \iota=p \circ \iota=\id$.

\begin{df} (i) We say that a subset $S$ of $X$ is $\RP$-connected if
$S \cap b$ is connected or empty for each orbit $b$ of $\mu$.

(ii)  Let $S$ be a subset of $X$. We set $\RP S=\mu(S,\RP).$

(iii) Let $S$ be a subset of $X$. Then $S$ is conic if $S=\RP
S$. i.e. $S$ is invariant under the action of $\RP$.
\end{df}

If $U \in \op(X)$, then $\RP U \in \op(X)$ because $\mu$ is open.

\begin{df} A sheaf of $k$-modules $F$ on $X_{sa}$ is conic if the restriction morphism
$\Gamma(\RP U;F) \to \Gamma(U;F)$ is an isomorphism for each
$\RP$-connected $U \in \op^c(X_{sa})$ with $\RP U \in
\op(X_{sa})$.
\begin{itemize}
\item[(i)]We denote by $\mod_{\RP}(k_{X_{sa}})$ the full subcategory of
$\mod(k_{X_{sa}})$ consisting of conic sheaves.
\item[(ii)] We denote by $D^b_{\RP}(k_{X_{sa}})$,
the full subcategory of $D^b(k_{X_{sa}})$ consisting of objects $F$
such that $H^jF$ belongs to $\mod_{\RP}(k_{X_{sa}})$ for all $j
\in \Z$.
\end{itemize}
\end{df}

\begin{oss}\label{noequivariant} Let $X$ be a real analytic manifold endowed with a subanalytic
action $\mu$ of $\RP$.
As in classical sheaf theory
one can define the subcategory $\mod^\mu(k_{X_{sa}})$ of
$\mod(k_{X_{sa}})$ consisting of sheaves satisfying $\imin \mu F
\simeq \imin p F$. The categories $\mod^\mu(k_{X_{sa}})$ and
$\mod_{\RP}(k_{X_{sa}})$ are not equivalent in general.

In fact let $X=\R$, set $X^+=\{x \in \R;\, x>0\}$ and let $\mu$ be
the natural action of $\RP$ (i.e. $\mu(x,t)=tx$). The sheaf
$\rho_!k_{X^+} \simeq \lind {n \in \N}\rho_*k_{({1\over n},n)}$ belongs to $\mod^\mu(k_{X_{sa}})$ but not to $\mod_{\RP}(k_{X_{sa}})$. Indeed, it is easy to check that $\Gamma(X^+;\rho_!k_{X^+})=0$ while $\Gamma((a,b);\rho_!k_{X^+}) \simeq k$, $0<a<b$. Moreover, since $k_{X^+}$ is conic and the functors $\imin \mu$ and $\imin p$ commute with $\rho_!$ we have $\rho_!k_{X^+} \in \mod^\mu(k_{X_{sa}})$.
\end{oss}

\begin{df}We denote by $\op(X_{\RP})$ (resp. $\op(X_{sa,\RP})$) the full subcategory of $\op(X)$ consisting of conic (resp. conic
subanalytic) subsets, i.e. $U \in \op(X_{\RP})$ (resp. $U \in
\op(X_{sa,\RP})$) if $U \in \op(X)$ (resp. $U \in \op(X_{sa})$)
and it is invariant by the action of $\RP$.

We denote by $X_{\RP}$ (resp. $X_{sa,\RP}$) the category
$\op(X_{\RP})$ (resp. $\op(X_{sa,\RP})$) endowed with the topology
induced by $X$ (resp. $X_{sa}$).
\end{df}

Let  $\eta: X \to X_{\RP}$, $\eta_{sa}: X_{sa} \to X_{sa,\RP}$ and $\rho_{\RP}:X_{\RP} \to X_{sa,\RP}$
be the natural morphisms of sites. We have a commutative diagram of sites
\begin{equation}\label{etarhoRP}
\xymatrix{X \ar[r]^{\rho} \ar[d]^{\eta} & X_{sa} \ar[d]^{\eta_{sa}} \\
X_{\RP} \ar[r]^{\rho_{\RP}} & X_{sa,\RP}.}
\end{equation}

We need to introduce the subcategory of coherent conic sheaves.

\begin{df}\label{2.2.2} Let $U \in \op(X_{\RP})$. Then $U$ is said to be
relatively quasi-compact if, for any covering $\{U_i\}_{i \in I}$
of $X_{\RP}$, there exists $J \subset I$ finite such that $U
\subset \bigcup_{i \in J} U_i.$ We write $U \subset\subset X_{\RP}$.

We  will denote by $\op^c(X_{\RP})$ the subcategory of
$\op(X_{\RP})$ consisting of relatively quasi-compact open
subsets.
\end{df}

 One can check easily that if $U \in \op^c(X)$, then $\RP
U \in \op^c(X_{\RP})$.



\begin{df} Let $F \in \mod(k_{X_{\RP}})$ and consider the family $\op(X_{sa,\RP})$.
\begin{itemize}
\item[(i)] $F$ is $X_{sa,\RP}$-finite if there exists an epimorphism
$G \twoheadrightarrow F$, with $G \simeq \oplus_{i \in I}
k_{U_i}$,
 $I$ finite and $U_i \in \op^c(X_{sa,\RP})$.
\item[(ii)] $F$ is $X_{sa,\RP}$-pseudo-coherent if for any morphism
$\psi:G \to F$, where $G$ is $X_{sa,\RP}$-finite, $\ker \psi$ is
$X_{sa,\RP}$-finite.
\item[(iii)] $F$ is $X_{sa,\RP}$-coherent if it is both $X_{sa,\RP}$-finite
and $X_{sa,\RP}$-pseudo-coherent.
\end{itemize}
We will denote by $\coh(X_{sa,\RP})$ the subcategory of
$\mod(k_{X_{\RP}})$ consisting of $X_{sa,\RP}$-coherent objects.
\end{df}

Replacing $\op^c(X_{sa})$ with $\op^c(X_{sa,\RP})$, we can adapt the
results of \cite{KS01,Pr1} and we get the following result
(see \cite{EP} for a detailed proof).

\begin{teo}\label{RPind} (i) Let $G \in \coh(X_{sa,\RP})$ and let
$\{F_i\}$ be a filtrant inductive system in $\mod(k_{X_{sa,\RP}})$.
Then we have an isomorphism
$$\lind i \Ho_{k_{X_{sa,\RP}}}(\rho_{\RP *}G,F_i) \iso
\Ho_{k_{X_{sa,\RP}}}(\rho_{\RP *}G,\lind i F_i).$$
Moreover the functor of direct image $\rho_{\RP*}$ associated to
the morphism $\rho_{\RP}$ in \eqref{etarhoRP} is fully faithful
and exact on $\coh(X_{sa,\RP})$.

(ii) Let $F \in \mod(k_{X_{sa,\RP}})$. There exists a small filtrant
inductive system $\{F_i\}_{i \in I}$ in $\coh(X_{sa,\RP})$ such
that $F \simeq \lind i \rho_{\RP *}F_i$.
\end{teo}

\begin{nt} Since $\rho_{\RP *}$ is fully faithful and exact
on $\coh(X_{sa,\RP})$, we can identify
$\coh(X_{sa,\RP})$ with its image in $\mod(k_{X_{sa,\RP}})$. When there
is no risk of confusion we will write $F$ instead of $\rho_{\RP *}F$,
for $F \in \coh(X_{sa,\RP})$.\\
\end{nt}

Let us consider the category $\mod_{\RP}(k_{X_{sa}})$ of conic
sheaves on $X_{sa}$. The restriction of $\eta_{sa*}$ induces a
functor denoted by $\widetilde{\eta}_{sa*}$ and we obtain a
diagram
\begin{equation}\label{sarho}
\xymatrix{\mod_{\RP}(k_{X_{sa}})  \ar[d] \ar@ <2pt>
[rr]^{\widetilde{\eta}_{sa*}} &&
  \mod(k_{X_{sa,\RP}}) \ar@ <2pt> [ll]^{\imin \eta_{sa}} \\
\mod(k_{X_{sa}}) \ar[urr]_{\eta_{sa*}} && }
\end{equation}

Now assume the hypothesis below:
\begin{equation}\label{hypsa}
  \begin{cases}
 \text{(i) every $U \in \op^c(X_{sa})$ has a finite covering consisting }\\
 \text{\ \ of $\RP$-connected subanalytic open subsets,}\\
 \text{(ii) for any $U \in \op^c(X_{sa})$ we have $\RP U \in \op(X_{sa})$,}\\
 \text{(iii) for any $x \in X$\ \ the set $\RP x$ is contractible,}\\
 \text{(iv) there exists a covering $\{V_n\}_{n \in \N}$ of $X_{sa}$ such that}\\
 \text{\ \ $V_n$ is $\RP$-connected and $V_n \subset\subset V_{n+1}$ for each $n$}.
  \end{cases}
\end{equation}

Let $U \in \op(X_{sa})$ such that $\RP U$ is still subanalytic.
Let $\varphi$ be the natural map from $\Gamma(\RP U;F)$ to
$\Gamma(U;\imin \eta_{sa} F)$ defined by
\begin{equation}\label{varphisa}
\begin{array}{ccc}
\Gamma(\RP U;F) & \to & \Gamma(\RP U; \eta_{sa*} \imin \eta_{sa} F) \\
  & \simeq & \Gamma(\RP U; \imin \eta_{sa} F) \\
 & \to & \Gamma(U; \imin \eta_{sa} F).
\end{array}
\end{equation}

\begin{prop}\label{*-1Usa} Let $F \in \mod(k_{X_{sa,\RP}})$.
Let $U \in \op(X_{sa})$, assume that $U$ is $\RP$-connected and that $\RP U$ is still subanalytic. Then the morphism
$\varphi$ defined by \eqref{varphisa} is an isomorphism.
\end{prop}
\dim\ \ (i) Assume that $U \in \op^c(X_{sa})$ is $\RP$-connected. Let $F \in \mod(k_{X_{sa,\RP}})$,
then $F=\lind i \rho_{\RP*} F_i$, with $F_i \in
\coh(X_{sa,\RP})$. We have the chain of
isomorphisms
\begin{eqnarray*}
\Ho_{k_{X_{sa}}}(k_U,\imin \eta_{sa}\lind i \rho_{\RP*}F_i) &
\simeq & \Ho_{k_{X_{sa}}}(k_U,
\lind i \rho_* \imin \eta F_i)\\
& \simeq & \lind i \Ho_{k_X}(k_U, \imin \eta F_i)\\
& \simeq & \lind i \Ho_{k_{X_{\RP}}}(k_{\RP U}, F_i)\\
& \simeq &  \Ho_{k_{X_{sa,\RP}}}(k_{\RP U},\lind i \rho_{\RP *}
F_i),
\end{eqnarray*}
where the first isomorphism follows since $\imin \eta_{sa} \circ
\rho_{\RP*} \simeq \rho_* \circ \imin \eta$ and the third one
follows from the equivalence between conic sheaves on $X$ and
sheaves on $X_{\RP}$. In the fourth isomorphism we used the fact
that $\RP U \in \op^c(X_{sa,\RP})$.

(ii) Let $U \in \op(X_{sa})$ be $\RP$-connected. Let $\{V_n\}_{n
\in \N} \in \cov(X_{sa})$ be a covering of $X$ as in \eqref{hypsa}
(iv) and set $U_n=U \cap V_n$. We have
\begin{equation}\label{ML0}
\Gamma(U;\imin \eta_{sa}F) \simeq \lpro n\Gamma(U_n;\imin \eta_{sa}F) \simeq \lpro n\Gamma(\RP U_n;F)
\simeq \Gamma(\RP U;F).
\end{equation} \qed

\begin{teo}\label{ex*sa} The functors $\widetilde{\eta}_{sa*}$ and $\imin \eta_{sa}$ in \eqref{sarho}
  are equivalences of ca\-te\-go\-ries inverse to each
  others.
\end{teo}
\dim\ \ (i) Let $F \in \mod_{\RP}(k_{X_{sa}})$, and let $U \in
\op^c(X_{sa})$ be $\RP$-connected. We have
$$\Gamma(U;F) \simeq  \Gamma(\RP U;F) \simeq
\Gamma(\RP U;\widetilde{\eta}_{sa*}F) \simeq \Gamma(U;\imin
\eta_{sa} \widetilde{\eta}_{sa*}F).$$ The third isomorphism
follows from Proposition \ref{*-1Usa}. Then \eqref{hypsa} (i)
implies $\imin \eta_{sa} \widetilde{\eta}_{sa*} \simeq \id$.

(ii) For any $U \in \op^c(X_{sa,\RP})$ we have:
$$\Gamma(U;\eta_{sa*}\imin \eta_{sa} F) \simeq \Gamma(U;\imin \eta_{sa} F) \simeq \Gamma(U;F)$$ where
the second isomorphisms follows from Proposition \ref{*-1Usa}. This implies $\eta_{sa*} \imin \eta_{sa} \simeq \id$.\\
\qed

\begin{nt} Since $\imin {\eta_{sa}}$ is fully faithful and exact
we will often identify $\coh(X_{sa,\RP})$ with its image in
$\mod_{\RP}(k_{X_{sa}})$. Hence, for $F \in \coh(X_{sa,\RP})$ we
shall often write $F$ instead of $\imin {\eta_{sa}}F$.
\end{nt}

Thanks to Theorem \ref{RPind} we can give another description of the category of conic sheaves.

\begin{teo}\label{indRP} Let $F \in \mod_{\RP}(k_{X_{sa}})$. Then there exists a small filtrant system $\{F_i\}$ in $\coh(X_{sa,\RP})$ such that  $F \simeq \lind i \rho_* \imin \eta F_i$.
\end{teo}

Assume \eqref{hypsa}. Injective and quasi-injective objects of $\mod(k_{X_{sa}})$ are not contained in $\mod_{\RP}(k_{sa})$. 
For this reason we are going to introduce a subcategory which is useful when we try to find acyclic resolutions.


\begin{lem} Assume that $X$ satisfies \eqref{hypsa}. Then the following property is satisfied:
\begin{equation}\label{Dhypsa}
  \begin{cases}
    \text{ each finite covering of an $\RP$-connected $U \in \op^c(X_{sa})$ }\\
    \text{ has a finite refinement $\{V_i\}_{i=1}^n$ such that each ordered}\\
    \text{ union $\bigcup_{i=1}^jV_i$ is $\RP$-connected for each $j \in \{1,\ldots,n\}$}.
  \end{cases}
\end{equation}
\end{lem}
\dim\ \ Let $U \in \op^c(X_{sa})$ be $\RP$-connected. Then each finite covering of $U$ admits a finite refinement consisting of $\RP$-connected open subanalytic subsets. Let $\{U_i\}_{i=1}^n$ be a finite covering of $U$, $U_i \in \op^c(X_{sa})$ $\RP$-connected for each $i$. We will construct a refinement satisfying \eqref{Dhypsa}.

For $k=1,\dots,n$ set $V_{k11}:=U_k$ and $V_{k1i}:=U_{\sigma(i)} \cap \RP(U_k \cap U_{\sigma(i)})$ for $i=2,\dots,n$ and $\sigma(i)=i-1$ if $i\leq k$, $\sigma(i)=i$ if $i>k$. Then set $U_{k2}:=\bigcup_{i=1}^nV_{k1i}$ and $V_{k2i}:=U_{\sigma(i)} \cap \RP(U_{k2} \cap U_{\sigma(i)})$. For $j=1,\dots, n$ define recursively $U_{kj}=\bigcup_{\ell=1}^j\bigcup_{i=1}^nV_{k\ell i}$ and $V_{kji}=U_{\sigma(i)} \cap \RP(U_{kj} \cap U_{\sigma(i)})$. Remark that 
$\bigcup_{p=1}^j\bigcup_{\ell=1}^n\bigcup_{i=1}^nV_{p\ell i}= \bigcup_{p=1}^j\RP U_p \cap U$. By Lemma \ref{RPconn} below all the sets $V_{kji}$ are $\RP$-connected and $\{V_{kji}\}_{i,k,j}$ is a refinement of $\{U_i\}_i$ satisfying \eqref{Dhypsa} (with the lexicographic order). \\
\qed

\begin{lem}\label{RPconn} Assume that $X$ satisfies \eqref{hypsa} (iii). Let $U,V,W$ be open and $\RP$-connected. Then $U \cup (V \cap \RP(U \cap V)) \cup (W \cap \RP(U \cap W))$ is $\RP$-connected.
\end{lem}
\dim\ \ In what follows, when we write $\RP x$ we suppose that $\RP x \simeq \R$. If $\RP x = x$ everything becomes obvious.

(i) First remark that $U \cap V$ (resp. $U \cap W$, $V \cap W$) is $\RP$-connected. Indeed, let $x_1 \in U \cap \RP x$, $x_2 \in V \cap \RP x$ for some $x \in X$. Then $x_1=\mu(x,a)$, $x_2=\mu(x,b)$.  Every path in $\RP x$ connecting $x_1$ and $x_2$ contains $\mu(x,[a,b])$. Since $U$ and $V$ are $\RP$-connected then $U \cap V \supset \mu(x,[a,b])$.

(ii) Now let us prove that $U \cup (V \cap \RP(U \cap V))$ is $\RP$-connected. Let $x_1,x_2 \in U \cup (V \cap \RP(U \cap V)) \cap \RP x$ for some $x \in X$. Then $x_1=\mu(x,a)$, $x_2=\mu(x,b)$. We want to prove that $\mu(x,[a,b]) \subset U \cup (V \cap \RP(U \cap V))$. If $x_1,x_2 \in U$ it follows since $U$ is $\RP$-connected and if $x_1,x_2 \in V \cap \RP(U \cap V)$ it follows from (i). So we may assume that $x_1 \in U$ and $x_2 \in V \cap \RP(U \cap V)$. Since $U$ is $\RP$-connected and $x_2 \in \RP x_1$, there exists $y=\mu(x,c) \in U \cap V$. Then $\mu(x,[a,c]) \subset U$. In the same way $\mu(x,[b,c]) \subset V \cap \RP(U \cap V)$ and hence $\mu(x,[a,c] \cup [b,c]) \subset U \cup (V \cap \RP(U \cap V))$.

(iii) Let us show that $U \cup (V \cap \RP(U \cap V)) \cup (W \cap \RP(U \cap W))$ is $\RP$-connected. Let $x_1,x_2 \in U \cup (V \cap \RP(U \cap V)) \cup (W \cap \RP(U \cap W)) \cap \RP x$ for some $x \in X$. Then $x_1=\mu(x,a)$, $x_2=\mu(x,b)$. We want to prove that $\mu(x,[a,b]) \subset U \cup (V \cap \RP(U \cap V)) \cup (W \cap \RP(U \cap W))$. By (i) and (ii) we may reduce to the case $x_1 \in V$, $x_2 \in W$. As in (ii), there exist $y_1=\mu(x,c) \in U \cap V$ and $y_2=\mu(x,d) \in U \cap W$. Then $\mu(x,[c,d]) \in U$, $\mu(x,[a,c]) \subset V \cap \RP(U \cap V)$ and $\mu(x,[b,d]) \subset W \cap \RP(U \cap W)$. Hence $\mu(x,[c,d] \cup [a,c] \cup [b,d]) \in U \cup (V \cap \RP(U \cap V)) \cup (W \cap \RP(U \cap W))$ and the result follows.\\
\qed

\begin{df} A sheaf $F \in \mod(k_{X_{sa}})$ is $\RP$-quasi-injective if for each $\RP$-connected $U \in \op^c(X_{sa})$ the restriction morphism $\Gamma(X;F) \to \Gamma(U;F)$ is surjective.
\end{df}


Remark that the functor $\imin {\eta_{sa}}$ sends injective objects of $\mod(k_{X_{sa,\RP}})$ to $\RP$-quasi-injective objects since $\Gamma(U;\imin {\eta_{sa}}F) \simeq \Gamma(\RP U;F)$ if $U \in \op^c(X_{sa})$ is $\RP$-connected.
Moreover the category of $\RP$-injective objects is cogenerating since injective objects are cogenerating in $\mod(k_{X_{sa}})$. Once we have \eqref{Dhypsa} and \eqref{hypsa} (iv) it is easy to prove Propositions \ref{conflU} and \ref{3rp} below in the same way as the corresponding classical results for c-soft sheaves (see Propositions 2.5.8, 2.5.10 and Corollary 2.5.9 of \cite{KS90}).

\begin{prop} \label{conflU} Let  $\exs{F'}{F}{F''}$ be an exact sequence in
$\mod(k_{X_{sa}})$ and assume that $F'$ is $\RP$-quasi-injective. Let
$U \in \op(X_{sa})$ be $\RP$-connected. Then the sequence
$$\exs{\Gamma(U;F')}{\Gamma(U;F)}{\Gamma(U;F'')}$$
is exact.
\end{prop}

\begin{prop}\label{3rp} Let $F',F$ be $\RP$-quasi-injective and consider the exact
sequence  $\exs{F'}{F}{F''}$ in $\mod(k_{X_{sa}})$. Then $F''$ is
$\RP$-quasi-injective.
\end{prop}

It follows from the preceding results that

\begin{prop} \label{RPqinjinj} $\RP$-quasi-injective objects are injective with respect to the functor $\Gamma(U;\cdot)$, with $U \in \op(X_{sa})$ and $\RP$-connected.
\end{prop}

\begin{cor} $\RP$-quasi-injective objects are $\eta_{sa*}$-injective.
\end{cor}

\begin{teo}\label{Rex*sa} The categories $D^b(k_{X_{sa,\RP}})$ and
$D^b_{\RP}(k_{X_{sa}})$ are equivalent.
\end{teo}
\dim\ \ In order to prove this statement, it is enough to show
that $\imin \eta_{sa}$ is fully faithful. Let $F \in D^b(k_{X_{sa,\RP}})$ and let $F'$ be an injective complex quasi-isomorphic to $F$. Since $\imin {\eta_{sa}}$ sends injective objects to $\RP$-quasi-injective objects which are $\eta_{sa*}$-injective, we have
$R\eta_{sa*} \imin {\eta_{sa}} F \simeq \eta_{sa*} \imin {\eta_{sa}} F'
\simeq F' \simeq  F$.
This implies $R\eta_{sa*} \imin \eta_{sa} \simeq \id$, hence $\imin {\eta_{sa}}$ is fully faithful.\\
\qed

\begin{cor}\label{preceding} Let $F \in D^b_{\RP}(k_{X_{sa}})$ and let $U \in \op(X_{sa})$ be $\RP$-connected. Then ${\rm R}\Gamma(\RP U;F) \iso {\rm R}\Gamma(U;F)$.
\end{cor}

Hence for each
$F \in D^b_{\RP}(k_{X_{sa}})$ we have  $F \simeq \imin
{\eta_{sa}}F'$ with $F' \in D^b(k_{X_{sa,\RP}})$.



\begin{oss}
Thanks to these results, in order to prove that a morphism $F \to G$ in $D^b_{\RP}(k_{X_{sa}})$ is an isomorphism, it is enough to check that ${\rm R}\Gamma(U;F) \iso {\rm R}\Gamma(U;G)$ for each $U \in \op(X_{sa,\RP})$.
\end{oss}

\begin{oss} It is easy to check that the six Grothendieck operations, except the functor of proper direct image, preserve conic subanalytic sheaves. We refer to \cite{Pr2} for a detailed exposition.
\end{oss}

\subsection{Conic sheaves on vector bundles}

Let $E \stackrel{\tau}{\to} Z$ be a real vector bundle, with
dimension $n$ over a real analytic manifold $Z$. Then $\RP$ acts
naturally on $E$ by multiplication on the fibers. We identify $Z$
with the zero-section of $E$ and denote by $i:Z \hookrightarrow E$
the embedding. We set $\pE=E \setminus Z$ and $\pittau:\pE \to
Z$ denotes the projection.

\begin{lem}\label{bierstone} The category $\op(E_{sa})$ satisfies
\eqref{hypsa}.
\end{lem}
\dim\ \ Let us prove \eqref{hypsa} (i). Let $\{V_i\}_{i \in \N}$
be a locally finite covering of $Z$ with $V_i \in \op^c(Z_{sa})$
such that $\imin \tau (V_i) \simeq \R^m \times \R^n$ and let
$\{U_i\}$ be a refinement of $\{V_i\}$ with $U_i \in
\op^c(Z_{sa})$ and $\overline{U_i} \subset V_i$ for each $i$. Then
$U$ is covered by a finite number of $\imin \tau (U_i)$ and $U
\cap \imin \tau(U_i)$ is relatively compact in $\imin \tau (V_i)$
for each $i$. We may reduce to the case $E \simeq \R^m \times
\R^n$.
 Let us consider the morphism
of manifolds
\begin{eqnarray*}
\varphi : \R^m \times \mathbb{S}^{n-1} \times \R & \to & \R^m \times \R^n\\
(z,\vartheta,r) & \mapsto & (z,r i(\vartheta)),
\end{eqnarray*}
where $i:\mathbb{S}^{n-1} \hookrightarrow \R^n$ denotes the
embedding. Then $\varphi$ is proper and subanalytic. The subset
$\imin \varphi(U)$ is subanalytic and relatively compact in $\R^m
\times  \mathbb{S}^{n-1} \times \R$. \\


(a) By a result of \cite{Wi05}, $\imin\varphi(U \setminus Z)$ admits a finite cover $\{W_j\}_{j \in J}$ such that the intersections of each $W_j$ with the fibers of $\pi:\R^m \times \mathbb{S}^{n-1} \times \R \to \R^m \times \mathbb{S}^{n-1}$ are contractible or empty. Then $\varphi(W_j)$ is an open subanalytic relatively compact $\RP$-connected subset of $\R^m \times \R^n$ for each $j$. In this way we obtain a finite covering of $U \setminus Z$ consisting of $\RP$-connected subanalytic open subsets.



(b) Let $p \in \pi(\imin\varphi(U \cap Z))$. Then $\imin \pi(p) \cap U$ is a disjoint union of intervals. Let us consider the interval $(m(p),M(p))$, $m(p)<M(p) \in \R$ containing $0$. Set $W_Z=\{(p,r) \in U ;\; m(p)<r<M(p)\}$. 
The set $W_Z$
is open subanalytic (it is a consequence of Proposition 1.2, Chapter 6 of \cite{VD98}), contains $\imin\varphi(U \cap Z)$ and its intersections with the fibers of $\pi$ are contractible. Then $\varphi(W_Z)$ is an open $\RP$-connected subanalytic neighborhood of $U \cap Z$ and it is contained in $U$.

By (a) there exists a finite covering $\{\varphi(W_j)\}_{j \in J}$ of $U
\setminus Z$ consisting of $\RP$-connected subanalytic open
subsets, and $\varphi(W_Z) \cup \bigcup_{j \in J}\varphi(W_j)=U$.\\

By Proposition 8.3.8 of \cite{KS90} the category $\op(E_{sa})$
also satisfies \eqref{hypsa} (ii). Moreover \eqref{hypsa} (iii)
and (iv) are clearly satisfied. \nopagebreak \newline \qed

Now let us consider $E$ endowed with the conic topology. In this situation, an object $U \in \op(E_{\RP})$ is the union of $\pU \in \op(\pE_{\RP})$ and $U_Z \in \op(Z)$ such that $\imin \pittau (U_Z) \subset \pU$. If $U,V \in \op(E_{\RP})$, then $U \subset\subset V$ if $U_Z \subset\subset V_Z$ in $Z$ and $\pU \subset\subset \pV$ in $\pE_{\RP}$ (this means that $\pi(\pU) \subset\subset \pi(\pV)$ in $\pE/\RP$, where $\pi:\pE \to \pE/\RP$ denotes the projection). \\


Applying Theorem \ref{Rex*sa} we have the following

\begin{teo} The categories $D^b_{\RP}(k_{E_{sa}})$ and $D^b(k_{E_{sa,\RP}})$ are
equi\-valent.
\end{teo}

Consider the subcategory $\mod^{cb}_{\rc,\RP}(k_E)$ of
$\mod_{\rc,\RP}(k_E)$ consisting of sheaves whose support is
compact on the base (i.e. $\tau(\supp(F))$ is compact in $Z$). The
restriction of $\imin \eta$ to $\coh(E_{sa,\RP})$ gives rise (see \cite{Pr2} for more details) to
an equivalence of categories
\begin{equation}\label{rhosa}
\imin {\overline{\eta}}:\coh(E_{sa,\RP}) \iso
\mod^{cb}_{\rc,\RP}(k_E).
\end{equation}

As a consequence of Theorem \ref{indRP} one has the following

\begin{teo} Let $F \in \mod_{\RP}(k_{E_{sa}})$. Then there exists a small filtrant system $\{F_i\}$ in $\mod_{\rc,\RP}^{cb}(k_E)$ such that  $F \simeq \lind i \rho_*F_i$.
\end{teo}



We end this section with this result, which will be
useful in the next section.

\begin{lem}\label{itausa} Let $F \in D^b_{\RP}(k_{E_{sa}})$. Then:
\begin{itemize}
\item[(i)] $R\tau_*F \simeq \imin i F$.
\item[(ii)] $R\tau_{!!}F \simeq i^!F$.
\end{itemize}
\end{lem}
\dim\ \ (i) The adjunction morphism defines $R\tau_*F \simeq \imin
i \imin \tau R\tau_* F \to \imin i F$. Let $V \in \op^c(Z_{sa})$.
Then
$$\lind {U \supset V} R^k\Gamma(U;F) \simeq
\lind {\substack{U \supset V \\ \tau(U)=V}}R^k\Gamma(U;F) \simeq
R^k\Gamma(\imin \tau(V);F) \simeq R^k\Gamma(V;R\tau_*F),$$ where
$U \in \op(E_{sa})$ and $\RP$-connected. The second isomorphism
follows from Corollary \ref{preceding}.

(ii) The adjunction morphism defines  $i^!F \to
i^!\tau^!R\tau_{!!}F \simeq R\tau_{!!}F$. Let $V \in
\op^c(Z_{sa})$, and let $K$ be a compact subanalytic
$\RP$-connected neighborhood of $V$ in $E$. Then $\imin \tau(V)
\setminus K$ is $\RP$-connected and subanalytic, and $\RP(\imin
\tau(V) \setminus K)=\imin \tau(V) \setminus Z$. By Corollary \ref{preceding}
we have the isomorphism $\mathrm{R}\Gamma(\imin
\tau(V);\mathrm{R}\Gamma_ZF) \simeq \mathrm{R}\Gamma(\imin \tau
(V);\mathrm{R}\Gamma_KF)$.

It follows from the definition of $R\tau_{!!}$ that for any $k \in
\Z$ and $V \in \op^c(Z_{sa})$ we have $R^k\Gamma(V;R\tau_{!!}F)
\simeq \lind K R^k\Gamma(\imin \tau (V);\mathrm{R}\Gamma_KF)$,
where $K$ ranges through the family of compact subanalytic
$\RP$-connected neighborhoods of $V$ in $E$.

On the other hand  for any $k \in \Z$ we have $R^k\Gamma(V;i^!F)
\simeq R^k\Ho(i_*k_V,F) \simeq R^k\Ho(i_*\imin i \imin \tau k_V,F)
\simeq R^k\Gamma(\imin \tau(V),\mathrm{R}\Gamma_ZF)$ and the
result follows. \nopagebreak \newline \qed




\subsection{Fourier-Sato transformation}\label{fouriersato}

Let $E \stackrel{\tau}{\to} Z$ be a real vector bundle, with
dimension $n$ over a real analytic manifold $Z$ and $E^*
\stackrel{\pi}{\to} Z$ its dual. We identify $Z$ as the
zero-section of $E$ and denote $i:Z \hookrightarrow E$ the
embedding, we define similarly $i:Z \hookrightarrow E^*$. We
denote by $p_1$ and $p_2$ the projections from $E \times_Z E^*$:
$$
\xymatrix{& & E \underset{Z}{\times} E^* \ar[lld]_{p_1} \ar[rrd]^{p_2} & &\\
E \ar[rrd]^{\tau} & & && E^* \ar[lld]_{\pi}\\
& & Z && &}
$$
We set
$$P:=\{(x,y) \in E \underset{Z}{\times} E^*;\; \langle x,y\rangle \geq 0\}$$
$$P':=\{(x,y) \in E \underset{Z}{\times} E^*;\; \langle x,y\rangle \leq 0\}$$
and we define the functors
$$\left\{
\begin{array}{l}
\Psi_{P'}=Rp_{1*}\circ \mathrm{R}\Gamma_{P'} \circ p_2^!:D^b_{\RP}(k_{E^*_{sa}}) \to D^b_{\RP}(k_{E_{sa}})\\
\\
\Phi_{P'}=Rp_{2!!}\circ (\cdot)_{P'} \circ p_1^{-1}:D^b_{\RP}(k_{E_{sa}}) \to D^b_{\RP}(k_{E^*_{sa}})\\
\\
\Psi_{P}=Rp_{2*}\circ \mathrm{R}\Gamma_{P} \circ p_1^{-1}:D^b_{\RP}(k_{E_{sa}}) \to D^b_{\RP}(k_{E^*_{sa}})\\
\\
\Phi_{P}=Rp_{1!!}\circ (\cdot)_{P} \circ
p_2^!:D^b_{\RP}(k_{E^*_{sa}}) \to D^b_{\RP}(k_{E_{sa}})
\end{array}
\right.$$

\begin{oss} These functors are well defined, more generally they
send subanalytic sheaves to conic subanalytic sheaves.
\end{oss}





\begin{lem}\label{suppsa} Let $F \in D^b_{\RP}(k_{E_{sa}})$. Then $\supp
  ((\mathrm{R}\Gamma_P(\imin {p_1}F))_{P'})$ is contained in $Z  \times_Z E^*$.
\end{lem}
\dim\ \ We may reduce to the case $F \in \mod_{\RP}(k_{E_{sa}})$.
Then $F=\lind i \rho_*F_i$, with $F_i \in
\mod^{cb}_{\rc,\RP}(k_E)$. We have
\begin{eqnarray*}
H^k(\mathrm{R}\Gamma_P (\imin {p_1}\lind i \rho_*F_i)_{P'}) &
\simeq & \lind i H^k(\mathrm{R}\Gamma_P(\imin {p_1} \rho_* F_i)_{P'}) \\
& \simeq & \lind i \rho_* H^k(\mathrm{R}\Gamma_P(\imin {p_1} F_i)_{P'}) \\
& \simeq & \lind i \rho_* (H^k(\mathrm{R}\Gamma_P (\imin {p_1}
F_i)_{P'}))_{Z  \times_Z E^*},
\end{eqnarray*}
where
the last isomorphism follows from Lemma 3.7.6 of \cite{KS90}.
\nopagebreak
\newline \qed

\begin{lem}\label{ABsa} Let $A$ and $B$ be two closed subanalytic subsets of $E$ such that $A \cup B=E$,
and let $F \in D^b(k_{E_{sa}})$. Then $\mathrm{R}\Gamma_A(F_B)
\simeq (\mathrm{R}\Gamma_AF)_B$.
\end{lem}

\dim\ \ We have a natural arrow $(\Gamma_AF)_B \to \Gamma_A(F_B)$,
and $R(\Gamma_AF)_B \simeq (\mathrm{R}\Gamma_AF)_B$ since
$(\cdot)_B$ is exact. Then we obtain a morphism
$(\mathrm{R}\Gamma_AF)_B \to \mathrm{R}\Gamma_A(F_B)$. It is
enough to prove that for any $k \in \Z$ and for any $F \in
\mod(k_{E_{sa}})$ we have $(R^k\Gamma_AF)_B \iso
R^k\Gamma_A(F_B)$. Since both sides commute with filtrant $\Lind$,
we may assume $F \in \mod^c_{\rc}(k_E)$. Then the result follows
from the corresponding one for classical sheaves. \nopagebreak
\newline \qed

\begin{prop} The two functors $\Phi_{P'},\Psi_P:D^b_{\RP}(k_{E_{sa}}) \to
  D^b_{\RP}(k_{E^*_{sa}})$ are isomorphic.
\end{prop}
\dim\ \ We have the chain of isomorphisms:
\begin{eqnarray*}
\Phi_{P'}  F & = & Rp_{2!!}(\imin {p_1}F)_{P'} \\
 & \simeq & Rp_{2!!}\mathrm{R}\Gamma_P((\imin {p_1}F)_{P'}) \\
 & \simeq & Rp_{2!!}(\mathrm{R}\Gamma_P(\imin {p_1}F))_{P'} \\
 & \simeq & Rp_{2*}(\mathrm{R}\Gamma_P(\imin {p_1}F))_{P'} \\
 & \simeq & Rp_{2*}\mathrm{R}\Gamma_P(\imin {p_1}F).
\end{eqnarray*}
The first isomorphism follows from Lemma \ref{itausa} (ii), the
second one from Lemma \ref{ABsa}, the third one from Lemma
\ref{suppsa} and the last one from Lemma \ref{itausa} (i).
\nopagebreak \newline \qed

\begin{df} Let $F \in D^b_{\RP}(k_{E_{sa}})$.
\begin{itemize}
\item[$(i)$] The Fourier-Sato transform is the functor $$(\cdot)^{\land}:D^b_{\RP}(k_{E_{sa}}) \to D^b_{\RP}(k_{E_{sa}^*})$$
$$F^{\land}=\Phi_{P'}F \simeq \Psi_PF.$$
\item[$(ii)$] The inverse Fourier-Sato transform is the functor $$(\cdot)^{\vee}:D^b_{\RP}(k_{E_{sa}^*}) \to D^b_{\RP}(k_{E_{sa}})$$
$$F^{\vee}=\Psi_{P'}F \simeq \Phi_PF.$$
\end{itemize}
\end{df}

It follows from definition that the functors ${}^\land$ and
${}^\vee$ commute with $R\rho_*$ and $\imin \rho$. We have
quasi-commutative diagrams
$$
\xymatrix{D^b_{\RP}(k_E) \ar[d]^{R\rho_*} \ar@ <2pt> [r]^{\land} & D^b_{\RP}(k_{E^*}) \ar@ <2pt> [l]^{\vee} \ar[d]^{R\rho_*} \\
D^b_{\RP}(k_{E_{sa}}) \ar@ <2pt> [r]^{\land} &
D^b_{\RP}(k_{E_{sa}^*}) \ar@ <2pt> [l]^{\vee}} \ \ \ \
\xymatrix{D^b_{\RP}(k_E) \ar@ <2pt> [r]^{\land} & D^b_{\RP}(k_{E^*}) \ar@ <2pt> [l]^{\vee} \\
D^b_{\RP}(k_{E_{sa}}) \ar@ <2pt> [r]^{\land} \ar[u]^{\imin \rho} &
D^b_{\RP}(k_{E_{sa}^*}) \ar@ <2pt> [l]^{\vee} \ar[u]^{\imin
\rho}.}
$$

This implies that these functors are the extension to conic
subanalytic sheaves of the classical Fourier-Sato and inverse
Fourier-Sato transforms.

\begin{teo}\label{fouriersa} The functors ${}^{\wedge}$ and ${}^{\vee}$ are equivalence of categories, inverse to each others. In particular we have
$$\Ho_{D^b_{\RP}(k_{E_{sa}})}(F,G) \simeq \Ho_{D^b_{\RP}(k_{E_{sa}^*})}(F^\wedge,G^\wedge).$$
\end{teo}
\dim\ \ Let $F \in D^b_{\RP}(k_{E_{sa}})$. The functors
${}^{\wedge}$ and ${}^{\vee}$ are adjoint functors, then we have a
morphism $F \to F^{\land \vee}$. To show that it induces an
isomorphism it is enough to check that $\mathrm{R}\Gamma(U;F) \to
\mathrm{R}\Gamma(U;F^{\land \vee})$ is an isomorphism on a basis
for the topology of $E_{sa}$. Hence we may assume that $U$ is
$\RP$-connected. By Proposition \ref{preceding} we may suppose that
$U$ is an open su\-ba\-na\-ly\-tic cone of $E$. we have the chain
of isomorphisms:
\begin{eqnarray*}
\Rh(k_U,F^{\land \vee}) & = & \Rh(k_U,\Psi_{P'}\Phi_{P'}F)\\
 & \simeq & \Rh(\Phi_{P'}k_U,\Phi_{P'}F)\\
 & \simeq & \Rh(\Phi_{P'}k_U,\Psi_{P}F)\\
 & \simeq & \Rh(\Phi_{P}\Phi_{P'}k_U,F)\\
 & \simeq & \Rh(k_U,F),
\end{eqnarray*}
where the last isomorphism follows from Theorem 3.7.9 of
\cite{KS90} and from the fact that the functors ${}^\land$ and
${}^\vee$ commute with $R\rho_*$. Similarly we can show that for
$G \in D^b_{\RP}(k_{E_{sa}^*})$ we have an isomorphism $G^{\vee
\land} \iso G$. \nopagebreak \newline \qed



\begin{oss} In the complex case we have the same result with
$$P:=\{(x,y) \in E \underset{Z}{\times} E^*, {\rm Re}\langle x,y\rangle \geq 0\},$$
$$P':=\{(x,y) \in E \underset{Z}{\times} E^*, {\rm Re}\langle x,y\rangle \leq 0\}.$$
\end{oss}

\begin{oss} The Fourier-Sato isomorphism can be extended to the case of ind-sheaves (see \cite{KS01} for  complete exposition). On $E_{\RP}$ one can define the category $\ind(k_{E_{\RP}})$ of conic ind-sheaves: $F \in \ind(k_{E_{\RP}})$ if it is a filtrant ind-limit
of $F_i \in \mod^{cb}(k_{E_{\RP}})$ (i.e. with compact support on the base). With slight modifications to the results of this section one can extend the Fourier-Sato transform to this setting and prove that it induces an equivalence between the categories $D^b(\ind(k_{E_{\RP}}))$ and $D^b(\ind(k_{E^*_{\RP}}))$.
\end{oss}

\section{Laplace transform}

As an application of the preceding constructions we introduce the
conic sheaves of tempered and Whitney holomorphic functions in order to give a sheaf theoretical interpretation of the Laplace isomorphism of \cite{KS97}. We
refer to \cite{Ka84, KS96} for the definition of the functors of
temperate and formal cohomology and to \cite{KS97} for the action
of the Laplace transform on temperate and formal cohomology.

\subsection{Review on temperate cohomology}

From now on, the base sheaf is $\CC$. Let $M$ be a real analytic
manifold. One denotes by $\db_M$ the sheaf of Schwartz's
distributions, and by $\D_M$ the sheaf of finite order
differential operators with analytic coefficients. In \cite{Ka84}
the author defined the functor
$$\tho(\cdot,\db_M):\mod_{\rc}(\CC_M) \to \mod(\D_M)$$
in the following way: let $U$ be a subanalytic subset of $M$ and
$Z=M \setminus U$. Then the sheaf $\tho(\CC_U,\db_M)$ is defined
by the exact sequence
$$\exs{\Gamma_Z\db_M}{\db_M}{\tho(\CC_U,\db_M)}.$$
This functor is exact and extends as functor in the derived
category, from $D^b_{\rc}(\CC_M)$ to $D^b(\D_M)$. Moreover the
sheaf $\tho(F,\db_M)$ is soft for any $\R$-constructible sheaf
$F$.\\

Let us denote by $\C^\infty_M$ the sheaf of $\C^\infty$-functions
and  let $Z$ be a closed subset of $M$. One denotes by
$\II^\infty_{M,Z}$ the sheaf of $\C^\infty$-functions on $M$
vanishing up to infinite order on $Z$.

\begin{df} A Whitney function on a closed subset $Z$ of $M$ is an indexed
family $F=(F^k)_{k\in \N^n}$ consisting of continuous functions on
$Z$ such that $\forall m \in \N$, $\forall k \in \N^n$, $|k| \leq
m$, $\forall x \in Z$, $\forall \varepsilon >0$ there exists a
neighborhood $U$ of $x$ such that $\forall y,z \in U \cap Z$

$$\left|F^k(z)-\sum_{|j+k|\leq m}{(z-y)^j \over
j!}F^{j+k}(y)\right| \leq \varepsilon d(y,z)^{m-|k|}.$$ We denote
by $W^\infty_{M,Z}$ the space of  Whitney $\C^\infty$-functions on
$Z$. We denote by $\W^\infty_{M,Z}$ the sheaf $U \mapsto
W_{U,U\cap Z}^\infty$.
\end{df}

In \cite{KS96} the authors defined the functor
$$\cdot \wtens \C^\infty_M:\mod_{\rc}(\CC_M) \to \mod(\D_M)$$
in the following way: let $U$ be a subanalytic open subset of $M$
and $Z=M \setminus U$. Then $\CC_U \wtens
\C^\infty_M=\II^\infty_{M,Z}$, and $\CC_Z \wtens
\C^\infty_M=\W^\infty_{M,Z}$. This functor is exact and extends as
a functor in the derived category, from $D^b_{\rc}(\CC_M)$ to
$D^b(\D_M)$. Moreover the sheaf $F \wtens \C^\infty_M$ is soft for
any $\R$-constructible sheaf $F$.\\

Now let $X$ be a complex manifold, $X^\R$ the underlying real
analytic manifold and $\overline{X}$ the complex conjugate
manifold. The product $X \times \overline{X}$ is a
complexification of $X^\R$ by the diagonal embedding $X^\R
\hookrightarrow X \times \overline{X}$. One denotes by $\OO_X$ the
sheaf of holomorphic functions and by $\D_X$ the sheaf of finite
order differential operators with holomorphic coefficients. For $F
\in D^b_{\rc}(\CC_X)$ one sets
\begin{eqnarray*}
\tho(F,\OO_X) & = & \rh_{\D_{\overline{X}}}(
\OO_{\overline{X}},\tho(F,\db_{X^\R})), \\
F \wtens \OO_X & = & \rh_{\D_{\overline{X}}}( \OO_{\overline{X}},F
\wtens \C^\infty_{X^\R}),
\end{eqnarray*}
and these functors are called the functors of temperate and formal
cohomology respectively.

\subsection{The Weyl algebra}

Let $E$ be  complex vector space of finite dimension $n$. Let us
denote by $\OO(E)$ the polynomial ring on $E$. We denote by $D(E)$
the Weyl algebra on $E$, that is, the ring of differential
operators with coefficient in $\OO(E)$.

The Fourier transform $\land:D(E) \to D(E^*)$ and the inverse
Fourier transform $\vee:D(E^*) \to D(E)$ induce isomorphisms which
are defined as follows: let $(z_1,\ldots,z_n)$ and
$(\zeta_1,\ldots,\zeta_n)$ be two systems of coordinates in $E$
and $E^*$ respectively. We have
$$z_i^\land=-\partial_{\zeta_i} \ \ \text{and} \ \
\partial_{z_i}^\land=\zeta_i.$$
On the other hand we have
$$-z_i=\partial_{\zeta_i}^\vee \ \ \text{and} \ \
\partial_{z_i}=\zeta_i^\vee.$$
Let us consider the subanalytic site $E_{sa}$. We denote by
$D(E_{sa})$ (resp. $\OO(E_{sa}$)) the constant sheaf on $E_{sa}$
associated to $D(E)$ (resp. $\OO(E)$). We denote by
$\mod(D(E_{sa}))$ the category of $D(E_{sa})$-modules.


\subsection{The sheaves $\dbt_M$ and $\CWM$}

Let $M$ be a real analytic manifold.

\begin{df} One denotes by $\dbt_M$ the presheaf of tempered distributions on
$M_{sa}$ defined as follows:
$$U \mapsto \Gamma(M;\db_M)/\Gamma_{M\setminus U}(M;\db_M).$$
\end{df}

 As a consequence of the \L ojasievicz's inequalities
\cite{Lo59}, for $U,V \in \op(M_{sa})$ the sequence
$$\exs{\dbt_M(U \cup
V)}{\dbt_M(U)\oplus\dbt_M(V)}{\dbt_M(U \cap V)}$$ is exact. For each $U \in \op(M_{sa})$ the restriction morphism
$\Gamma(M;\dbt_M) \to \Gamma(U;\dbt_M)$ is surjective and ${\rm R}\Gamma(U;\dbt_M)$ is concentrated in degree zero. Moreover $\dbt_M$ is exact on $\mod_{\rc}(\CC_M)$, i.e. it is a quasi-injective object of $\mod(\CC_{M_{sa}})$. We have the
following result (see \cite{KS01}, Proposition 7.2.6)

\begin{prop} For each $F \in \mod_{\rc}(\CC_M)$ one has the isomorphism
$$\Ho(F,\dbt_M) \simeq \Gamma(M;\tho(F,\db_M)).$$
\end{prop}

Now let $X$ be a complex manifold, $X^\R$ the underlying real
analytic manifold and $\overline{X}$ the complex conjugate
manifold. One defines the sheaf $\ot_X \in D^b(\CC_{X_{sa}})$ of tempered
holomorphic functions as follows:
$$
\ot_X = \rh_{\rho_!\D_X}(\rho_!\OO_{\overline{X}},\dbtxr).
$$

The relation with the functor of temperate cohomology
is given by the following result

\begin{prop} For each $F \in D^b_{\rc}(\CC_X)$ one has the
isomorphism
$$
\tho(F,\OO_X) \simeq \imin \rho \rh(F,\ot_X).
$$
\end{prop}

Let $M$ be a real analytic manifold. As usual we set $D'(\cdot)=\rh(\cdot,\CC_M)$. We consider a slight generalization of the sheaf of Whitney $\C^\infty$-functions of \cite{KS01}.

\begin{df} Let $F \in \mod_{\rc}(\CC_M)$ and let $U \in \op(M_{sa})$. We define the presheaf $\CW_{M|F}$ as follows:
$$U \mapsto \Gamma(M;H^0D'\CC_U \otimes F \wtens \C^\infty_M).$$
\end{df}
Let $U,V \in \op(M_{sa})$, and consider the exact sequence
$$\exs{\CC_{U\cap V}}{\CC_U \oplus \CC_V}{\CC_{U\cup V}},$$
applying the functor $\ho(\cdot,\CC_M)=H^0D'(\cdot)$ we obtain
$$\lexs{H^0D'\CC_{U\cap V}}{H^0D'\CC_U \oplus H^0D'\CC_V}{H^0D'\CC_{U\cup V}},$$
applying the exact functors $\cdot \otimes F$, $\cdot\wtens \C^\infty_M$ and taking global sections
we obtain
$$\lexs{\CW_{M|F}(U\cup V)}{\CW_{M|F}(U) \oplus \CW_{M|F}(V)}{\CW_{M|F}(U\cap V)}.$$
This implies that $\CW_{M|F}$ is a sheaf on $M_{sa}$. Moreover if $U \in \op(M_{sa})$ is l.c.t., the morphism $\Gamma(M;\CW_{M|F}) \to \Gamma(U;\CW_{M|F})$ is surjective and $\mathrm{R}\Gamma(U;\CW_{M|F})$ is concentrated in degree zero. Let $\exs{F}{G}{H}$ be an exact sequence in $\mod_{\rc}(\CC_M)$, we obtain an exact sequence in $\mod(\CC_{M_{sa}})$
\begin{equation}\label{exsFGH}
\exs{\CW_{M|F}}{\CW_{M|G}}{\CW_{M|H}}.
\end{equation}

We can easily extend the sheaf $\CW_{M|F}$ to the case of $F \in D^b_{\rc}(\CC_M)$, taking a finite resolution of $F$ consisting of locally finite sums $\oplus \CC_V$, with $V$ l.c.t. in $\op^c(M_{sa})$. In fact, the sheaves $\CW_{M|\oplus \CC_V}$ form a complex quasi-isomorphic to $\CW_{M|F}$ consisting of acyclic objects with respect to $\Gamma(U;\cdot)$, where $U$ is l.c.t. in $\op^c(M_{sa})$.

As in the case of Whitney $\C^\infty$-functions one can prove that, if $G \in D^b_{\rc}(\CC_M)$ one has
$$\imin \rho \rh(G,\CW_{F|M}) \simeq D'G \otimes F \wtens \C^\infty_M.$$


\begin{es} Setting $F=\CC_M$ we obtain the sheaf of Whitney $\C^\infty$-functions. Let $U$ be an open subanalytic subset of $M$. Then $\CW_{M|\CC_U}$ is the sheaf of Whitney $\C^\infty$-functions vanishing on $M \setminus U$ with all their derivatives.
\end{es}

\begin{nt} Let $Z$ be a locally closed subanalytic subset of $M$. We set for short $\CW_{M|Z}$ instead of $\CW_{M|\CC_Z}$.
\end{nt}

Let $X$ be a complex manifold and let $Z$ be a complex submanifold of $X$. Let $F \in D^b_{\rc}(\CC_X)$. We denote by $\OW_{X|F} \in D^b(\CC_{X_{sa}})$ the object defined as follows:
$$
\OW_{X|F} := \rh_{\rho_!\D_{\overline{X}}}(\rho_!\OO_{\overline{X}},\CW_{X_\R|F}).
 $$
Let $\exs{F}{G}{H}$ be an exact sequence in $\mod_{\rc}(\CC_X)$. Then the exact sequence \eqref{exsFGH} gives rise to the distinguished triangle
\begin{equation}\label{dtFGH}
\dt{\OW_{X|F}}{\OW_{X|G}}{\OW_{X|H}}.
\end{equation}

The relation with the functor of formal cohomology
is given by the following result

\begin{prop} For each $F,G \in D^b_{\rc}(\CC_X)$ one has the
isomorphism
$$
D'F \otimes G \wtens \OO_X \simeq  \imin \rho \rh(F,\OW_{X|G}).
$$
\end{prop}

\subsection{Direct images for proper smooth morphisms}

Let $X$ be a complex manifold of complex dimension ${\rm dim}X$. Let $\Omega^t_X$ denote the subanalytic sheaf of tempered holomorphic forms of degree ${\rm dim}X$. Let $f:X \to Y$ be a morphism of complex manifolds. Let $\ddxy=\OO_X \otimes_{\imin f\OO_Y}\imin f\D_Y$ be the transfer bimodule of $f$. The inverse image of a $\D_Y$-module $\M$ is defined by
$$
\imin{\underline{f}}\M=\ddxy \underset{\imin f \D_Y}{\ltens} \imin f \M.
$$
 Set $d=\mathrm{dim}X-\mathrm{dim}Y$. We recall the following isomorphisms of \cite{KS01} and \cite{Pr2}:
\begin{eqnarray}
f^!\Omega^t_Y & \simeq & \Omega^t_X \underset{\rho_!\D_X}{\ltens} \rho_!\ddxy[d], \label{iminvot} \\
f^!\OWY & \simeq & \rh_{\rho_!\D_X}(\rho_!\ddxy,\OWX)[2d]. \label{iminvow}
\end{eqnarray}
Isomorphisms \eqref{iminvotDmod} and \eqref{dirimotDmod} below has been already proven in \cite{KS01} in the framework of ind-sheaves and can be obtained thanks to the equivalence between ind-$\R$-constructible sheaves and subanalyitic sheaves. Here we give a direct proof.

\begin{prop} \label{iminDmod} Let $\M \in D^b(\D_Y)$. There are natural isomorphisms
\begin{eqnarray}
f^!(\Omega^t_Y \underset{\rho_!\D_Y}{\ltens} \rho_!\M) & \simeq & \Omega^t_X \underset{\rho_!\D_X}{\ltens} \rho_!\underline{f}^{-1}\M[d], \label{iminvotDmod} \\
f^!\rh_{\rho_!\D_Y}(\rho_!\M,\OWY) & \simeq & \rh_{\rho_!\D_X}(\rho_!\underline{f}^{-1}\M,\OWX)[2d]. \label{iminvowDmod}
\end{eqnarray}
\end{prop}
\dim\ \ (i) We have the chain of isomorphisms
\begin{eqnarray*}
f^!(\Omega^t_Y \underset{\rho_!\D_Y}{\ltens} \rho_!\M) & \simeq & f^!\Omega^t_Y \underset{\rho_!\imin f\D_Y}{\ltens} \rho_!\imin f\M \\
& \simeq & \Omega^t_X \underset{\rho_!\D_X}{\ltens} \rho_!\ddxy \underset{\rho_!\imin f\D_Y}{\ltens} \rho_!\imin f\M[d] \\
& \simeq & \Omega^t_X \underset{\rho_!\D_X}{\ltens} \rho_!\underline{f}^{-1}\M[d],
\end{eqnarray*}
where the first isomorphism follows from Proposition 2.4.7 of \cite{Pr1} and the second one from \eqref{iminvot}.

(ii) We have the chain of isomorphisms
\begin{eqnarray*}
\lefteqn{f^!\rh_{\rho_!\D_Y}(\rho_!\M,\OWY)} \\ & \simeq & \rh_{\rho_!\imin f\D_Y}(\rho_!\imin f\M,f^!\OWY) \\
& \simeq & \rh_{\rho_!\imin f\D_Y}(\rho_!\imin f\M,\rh_{\rho_!\D_X}(\rho_!\ddxy,\OWX))[2d] \\
& \simeq & \rh_{\rho_!\D_X}(\rho_!\underline{f}^{-1}\M,\OWX))[2d],
\end{eqnarray*}
where the first isomorphism follows from the dual projection formula and the second one from \eqref{iminvow}.\\
\qed

We say that an $\OO_X$-module ${\cal F}$ is quasi-good if for any relatively compact open subset $U$, ${\cal F}$ is an increasing sequence of coherent $\OO_X|_U$-submodules. A $\D_X$-module is called quasi-good if it is quasi-good as an $\OO_X$-module. We denote by $D^b_{q-good}(\D_X)$ the full triangulated subcategory of $D^b(\D_X)$ consisting of objects with quasi-good cohomology. The direct image of a $\D_X$-module $\NN$ is defined by
$$
\underline{f}{}_*\NN=Rf_*(\ddyx \underset{\D_X}{\ltens} \NN ).
$$

\begin{prop} \label{dirimDmod} Suppose that $f:X \to Y$ is smooth and proper and let $\NN \in D^b_{q-good}(\D_X)$. Then we have natural isomorphisms
\begin{eqnarray}
Rf_*(\Omega^t_X \underset{\rho_!\D_X}{\ltens} \rho_! \NN) & \simeq & \Omega^t_Y \underset{\rho_!\D_Y}{\ltens} \rho_!\underline{f}{}_*\NN \label{dirimotDmod} \\
Rf_*\rh_{\rho_!\D_X}(\rho_!\NN,\OWX)[d] & \simeq & \rh_{\rho_!\D_Y}(\rho_!\underline{f}{}_*\NN,\OWY). \label{dirimowDmod}
\end{eqnarray}
\end{prop}
\dim\ \ We shall first find the morphisms and then show that they are isomorphisms.

${\rm (i)}_a$ We have the morphisms
\begin{eqnarray*}
\imin f\Omega^t_Y \underset{\rho_!\imin f\D_Y}{\ltens} \rho_!\ddyx & \simeq & f^!\Omega^t_Y \underset{\rho_!\imin f\D_Y}{\ltens} \rho_!\ddyx[-2d] \\
& \simeq & \Omega^t_X \underset{\rho_!\D_X}{\ltens} \rho_!\ddxy \underset{\rho_!\imin f\D_Y}{\ltens} \rho_!\ddyx[-d] \\
& \simeq & \rh_{\rho_!\D_X}(\rho_!\ddyx,\Omega^t_X) \underset{\rho_!\imin f\D_Y}{\ltens} \rho_!\ddyx \\
& \to & \Omega^t_X,
\end{eqnarray*}
where the first isomorphism follows from Proposition 2.4.9 of \cite{Pr1}, the second one follows from \eqref{iminvot}, and the third one from Proposition 4.19 of \cite{Ka03}.

${\rm (i)}_b$ We have the morphisms
\begin{eqnarray*}
\Omega^t_Y \underset{\rho_!\D_Y}{\ltens} \rho_!\underline{f}{}_*\NN & \simeq & \Omega^t_Y \underset{\rho_!\D_Y}{\ltens} \rho_!Rf_*(\ddyx \underset{\D_X}{\ltens} \NN) \\
& \simeq & \Omega^t_Y \underset{\rho_!\D_Y}{\ltens} Rf_*(\rho_!\ddyx \underset{\rho_!\D_X}{\ltens} \rho_!\NN) \\
& \simeq & Rf_*(\imin f\Omega^t_Y \underset{\rho_!\imin f\D_Y}{\ltens} (\rho_!\ddyx \underset{\rho_!\D_X}{\ltens} \rho_!\NN)) \\
& \simeq & Rf_*((\imin f\Omega^t_Y \underset{\rho_!\imin f\D_Y}{\ltens} \rho_!\ddyx ) \underset{\rho_!\D_X}{\ltens} \rho_!\NN) \\
& \to & Rf_*(\Omega^t_X \underset{\rho_!\D_X}{\ltens} \rho_!\NN),
\end{eqnarray*}
where
the second isomorphism follows since $f$ is smooth and
the last arrow follows from ${\rm (i)}_a$.

${\rm (i)}_c$ Let us prove \eqref{dirimotDmod}. Let $U \in \op^c(Y_{sa})$.
We have the chain of isomorphisms
\begin{eqnarray*}
{\rm R}\Gamma(U;\Omega^t_Y \underset{\rho_!\D_Y}{\ltens} \rho_!\underline{f}{}_*\NN) & \simeq & {\rm R}\Gamma(Y;\imin \rho{\rm R}\Gamma_U(\Omega^t_Y \underset{\rho_!\D_Y}{\ltens} \rho_!\underline{f}{}_*\NN)) \\
& \simeq & {\rm R}\Gamma(Y;\tho(\CC_U,\Omega_Y ) \underset{\D_Y}{\ltens} \underline{f}{}_*\NN) \\
& \simeq & {\rm R}\Gamma(Y;Rf_*(\tho(\CC_{\imin f(U)},\Omega_X ) \underset{\D_X}{\ltens} \NN)) \\
& \simeq & {\rm R}\Gamma(X;\imin \rho{\rm R}\Gamma_{\imin f(U)}(\Omega^t_X  \underset{\rho_!\D_X}{\ltens}\rho_! \NN)) \\
& \simeq & {\rm R}\Gamma(\imin f(U) ;\Omega_X^t  \underset{\rho_!\D_X}{\ltens} \rho_!\NN) \\
& \simeq & {\rm R}\Gamma(U;Rf_*(\Omega^t_X  \underset{\rho_!\D_X}{\ltens} \rho_!\NN)),
\end{eqnarray*}
where the third isomorphism follows from Theorem 7.3 of \cite{KS96}.\\

${\rm (ii)}_a$  We have the isomorphisms
\begin{eqnarray*}
f^!\OWY & \simeq &
 \rh_{\rho_!\D_X}(\rho_!\ddxy,\OWX)[2d] \\
& \simeq & \rho_!\ddyx \underset{\rho_!\D_X}{\ltens} \OWX[d],
\end{eqnarray*}
where the first isomorphism follows from \eqref{iminvow}  and the second one from Proposition 4.19 of \cite{Ka03}.

${\rm (ii)}_b$  We have the morphisms
\begin{eqnarray*}
\lefteqn{Rf_*\rh_{\rho_!\D_X}(\rho_!\NN,\OWX)[d]} \\ & \to & Rf_*\rh_{\rho_!\imin f\D_Y}(\rho_!(\ddyx\underset{\D_X}{\ltens} \NN),\rho_!\ddyx \underset{\rho_!\D_X}{\ltens} \OWX)[d] \\
& \simeq & Rf_*\rh_{\rho_!\imin f\D_Y}(\rho_!(\ddyx\underset{\D_X}{\ltens} \NN),f^! \OWY) \\
& \simeq & \rh_{\rho_!\D_X}(Rf_*\rho_!(\ddyx\underset{\D_X}{\ltens} \NN),\OWY) \\
& \simeq & \rh_{\rho_!\D_X}(\rho_!Rf_*(\ddyx\underset{\D_X}{\ltens} \NN),\OWY) \\
& \simeq & \rh_{\rho_!\D_X}(\rho_!\underline{f}{}_*\NN,\OWY),
\end{eqnarray*}
where the first isomorphism follows from ${\rm (ii)}_a$ and the third one follows since $f$ is smooth.

${\rm (ii)}_c$ Let us prove \eqref{dirimowDmod}. Let $U \in \op^c(Y_{sa})$.
We have the chain of isomorphisms
\begin{eqnarray*}
{\rm R}\Gamma(U;\rh_{\rho_!\D_Y}(\rho_!\underline{f}{}_*\NN,\OWY) & \simeq & \Rh_{\rho_!\D_Y}(\rho_!\underline{f}{}_*\NN,{\rm R}\Gamma_U\OWY) \\
& \simeq & \Rh_{\D_Y}(\underline{f}{}_*\NN,\imin \rho{\rm R}\Gamma_U\OWY) \\
& \simeq & \Rh_{\D_Y}(\underline{f}{}_*\NN,D'\CC_U \wtens \OO_Y) \\
& \simeq & \Rh_{\D_X}(\NN,\imin fD'\CC_U \wtens \OO_X)[d] \\
& \simeq & \Rh_{\D_X}(\NN,D'\imin f\CC_U \wtens \OO_X)[d] \\
& \simeq & \Rh_{\D_X}(\NN,\imin \rho{\rm R}\Gamma_{\imin f(U)} \OW_X)[d] \\
& \simeq & \Rh_{\rho_!\D_X}(\rho_!\NN,{\rm R}\Gamma_{\imin f(U)}\OWX)[d] \\
&  \simeq & {\rm R}\Gamma(\imin f(U);\rh_{\rho_!\D_X}(\rho_!\NN,\OWX))[d] \\
& \simeq & {\rm R}\Gamma(U;Rf_*\rh_{\rho_!\D_X}(\rho_!\NN,\OWX)[d],
\end{eqnarray*}
where the fourth isomorphism follows from Theorem 7.3 of \cite{KS96} and the fifth one follows since $f$ is smooth.\\
\qed

\begin{oss} 

There are similar isomorphisms for $\OW_{Y|F}$, $F \in D^b_{\rc}(\CC_Y)$, namely (with the same hypothesis as above)
\begin{eqnarray}
f^!\rh_{\rho_!\D_Y}(\rho_!\M,\OW_{Y|F}) & \simeq & \rh_{\rho_!\D_X}(\rho_!\underline{f}^{-1}\M,\OW_{X|\imin fF})[2d], \label{iminvowDmodF}\\
\rh(\rho_!\underline{f}{}_*\NN,\OW_{Y|F}) & \simeq & Rf_*\rh_{\rho_!\D_X}(\NN,\OW_{X|\imin fF})[d]. \label{dirimowDmodF}
\end{eqnarray}
The proof is the same as the one for \eqref{iminvowDmod} and \eqref{dirimowDmod}. We only considered the case $F=\CC_X$ to lighten notations.
\end{oss}

\subsection{The sheaves $\dbt_{E_{\RP}}$ and $\CW_{E_{\RP}}$}
Let $E$ be a $n$-dimensional real vector space, and let $i:E
\hookrightarrow P$ its projective compactification. If there is no
risk of confusion we will identify $E$ with its image $i(E)$. Let
us denote by $j:P \to E_{\RP}$ the natural morphism of sites, i.e.
$j^t(U)=U$ for each $U \in \op(E_{\RP})$.

\begin{df}\label{dbterp} We denote by $\dbt_{E_{\RP}}$ the sheaf on
$E_{sa,\RP}$ defined as follows: $\dbt_{E_{\RP}}:=j_*\dbt_P$.
\end{df}
Then $\dbt_{E_{\RP}}$ is a sheaf on $E_{sa,\RP}$, that is, it
belongs to $\mod_{\RP}(\CC_{E_{sa}})$. Moreover $\dbt_{E_{\RP}}$
extends to an exact functor on $\mod_{\rc,\RP}(\CC_E)$, i.e.
$\Ho(\cdot,\dbt_{E_{\RP}})$ is exact on $\mod_{\rc,\RP}(\CC_E)$.
In fact let $F \in \mod_{\rc,\RP}(\CC_E)$. We have
$$\Ho(F,\dbt_{E_{\RP}}) \simeq \Ho(F,j_*\dbt_P) \simeq \Ho(\imin
j F, \dbt_P),$$ and $\dbt_P$ is exact on $\mod_{\rc}(\CC_P)$ since
it is quasi-injective.
Let us consider $\Th(F,\db_E):=\Gamma(P;\tho(\imin jF,\db_P))$.
\begin{oss} It is isomorphic to $\Th(F,\db_E):=\Gamma(P;\tho(i_!F,\db_P))$
introduced in \cite{KS97}. In fact one can check easily that $i_!F
\simeq \imin j F$.
\end{oss}
We have the chain of isomorphisms
\begin{eqnarray*}
\Ho(F,\dbt_{E_{\RP}}) & \simeq & \Ho(\imin j F, \dbt_P) \\
& \simeq & \Gamma(P;\tho(\imin j F, \db_P)) \\
& \simeq & \Th(F,\db_E).
\end{eqnarray*}
It follows that for any $U \in \op(E_{sa,\RP})$ we have
$$\Gamma(U;\dbt_{E_{\RP}}) \simeq \Gamma(P;\tho(\imin j \CC_U, \db_P)) \simeq \Th(\CC_U,\db_E).$$
Hence the sections of $\dbt_{E_{\RP}}(U)$ are tempered along the
boundary of $U$ and at infinity. Denote by $D(E)$ the Weyl algebra
on $E$.
Then $\dbt_{E_{\RP}}$ is a $D(E_{sa})$-module.\\



\begin{df}\label{dbterp} We denote by $\CW_{E_{\RP}}$ the sheaf on
$E_{sa,\RP}$ defined as follows: $\CW_{E_{\RP}}:=j_*\CW_{P|E}$.
\end{df}
Then $\CW_{E_{\RP}}$ is a sheaf on $E_{sa,\RP}$, that is, it
belongs to $\mod_{\RP}(\CC_{E_{sa}})$. Moreover if $V$ is an open subanalytic l.c.t. cone then ${\rm R}\Gamma(V;\CW_{E_{\RP}})$
is concentrated in degree zero.
In fact we have
$$\Ho(\CC_V,\CW_{E_{\RP}}) \simeq \Ho(\CC_V,j_*\CW_{P|E}) \simeq \Ho(\CC_V, \CW_{P|E}),$$ and $\CW_{P|E}$ is $\Gamma(V;\cdot)$-acyclic.
Let $F \in D^b_{\rc,\RP}(\CC_E)$ and consider the functor
$F\stackrel{\rm W}{\otimes}\C^\infty_E:= {\rm R}\Gamma(P;\imin j F
\wtens \C^\infty_P)$.
\begin{oss} It is isomorphic to $F\stackrel{\rm W}{\otimes}\C^\infty_E:= {\rm R}\Gamma(P;i_! F
\wtens \C^\infty_P)$
introduced in \cite{KS97}.
\end{oss}
We have the chain of isomorphisms
\begin{eqnarray*}
\Ho(F,\CW_{E_{\RP}}) & \simeq & \Ho(\imin j F, \CW_{P|E}) \\
& \simeq & \Gamma(P;D'(\imin jF) \otimes \CC_E \wtens \C^\infty_P) \\
& \simeq & D'F\stackrel{\rm W}{\otimes}\C^\infty_E.
\end{eqnarray*}
It follows that for any $U \in \op(E_{sa,\RP})$ l.c.t. we have
$$\Gamma(U;\CW_{E_{\RP}}) \simeq \Gamma(U;\CW_{P|E}) \simeq  \Gamma(P;\CC_{\overline{U} \cap E} \wtens \C^\infty_P) \simeq \CC_{\overline{U}}\stackrel{\rm W}{\otimes}\C^\infty_E.$$
Hence the sections of $\CW_{E_{\RP}}(U)$ are Whitney functions $U$ rapidly decay at infinity. Denote by $D(E)$ the Weyl algebra
on $E$.
Then $\CW_{E_{\RP}}$ is a $D(E_{sa})$-module.\\


\subsection{Laplace transform}

Let $E$ be a $n$-dimensional complex vector space,   let $\overline{E}$ be the complex conjugate of $E$.

\begin{df} We define the conic sheaves of tempered and Whitney
holomorphic functions by the complexes
\begin{eqnarray*}
\ot_{E_{\RP}} & = &
\rh_{D(\overline{E}_{sa})}(\OO(\overline{E}_{sa}),\dbt_{E_{\RP}}),
\\
\OW_{E_{\RP}} & = &
\rh_{D(\overline{E}_{sa})}(\OO(\overline{E}_{sa}),\CW_{E_{\RP}}).
\end{eqnarray*}
\end{df}
It follows from the definition that for any $F \in
D^b_{\rc,\RP}(\CC_E)$
\begin{eqnarray*}
\Rh(F,\ot_{E_{\RP}}) & \simeq & \Th(F,\OO_E), \\
\Rh(F,\OW_{E_{\RP}}) & \simeq & D'F\overset{\rm W}{\otimes}\OO_E.
\end{eqnarray*}

Let $E$ be a $n$-dimensional complex vector space,   let $E^*$
the dual vector space. Let $i_1:E
\hookrightarrow P$  and $i_2:E^*
\hookrightarrow P^*$ denote the projective compactifications of $E$ and $E^*$. Set $i:E \times E \hookrightarrow P \times P^*$. Let $j_1:P
\hookrightarrow E_{\RP}$  and $j_2:P^*
\hookrightarrow E_{\RP}^*$ denote the natural morphisms of sites. We still denote by $p_1,p_2$ the projections from $P \times P^*$ to $P$ and $P^*$ respectively. We also still denote by ${}^\wedge$ and ${}^\vee$ the Fourier-Sato transforms on $P$ and $P^*$ with kernels $\CC_{i(P')}$ and $\CC_{i(P)}$ respectively.

\begin{lem} \label{compact} Let $G \in D^b(\CC_{P^*_{sa}})$. Then
$$
(Rj_{2*}G)^\vee \simeq Rj_{1*}(G^\vee).
$$
\end{lem}
\dim\ \ (i) Let $F \in D^b_{\rc,\RP}(\CC_E)$. Then we have the isomorphism \cite{KS97}
\begin{equation}\label{isocompact}
i_{2!}(F^\land) \simeq (i_{1!}F)^\land.
\end{equation}
Then, keeping in mind that $i_{k!}F \simeq \imin {j_k}F$, $k=1,2$ we have
\begin{eqnarray*}
\Rh(F,(j_{2*}G)^\vee) & \simeq & \Rh(i_{2!}(F^\land),G) \\
& \simeq & \Rh((i_{1!}F)^\land,G) \\
& \simeq & \Rh(F,j_{1*}(G^\vee)),
\end{eqnarray*}
where the first and second isomorphism follow by adjunction and the second one follows from \eqref{isocompact}.\\
 \qed

Set $S=P \times P^* \setminus i(E \times E^*)$ and let $\OO_{P \times P^*}(*S)$ be the sheaf of meromorphic functions whose poles are contained in $S$. Let us consider the $\D_{P \times P^*}$-modules
\begin{eqnarray*}
\LL & = & (\D_{P \times P^*}e^{-\langle x,y\rangle}) \underset{\OO_{P \times P^*}}{\ltens} \OO_{P \times P^*}(*S), \\
\LL' & = & (\D_{P \times P^*}e^{\langle x,y\rangle}) \underset{\OO_{P \times P^*}}{\ltens} \OO_{P \times P^*}(*S).
\end{eqnarray*}
$\LL$ is an holonomic $\D_{P \times P^*}$-module satisfying $\LL \simeq \LL \underset{\OO_{P \times P^*}}{\ltens} \OO_{P \times P^*}(*S).$ As an $\OO_{P \times P^*}(*S)$-module it is invertible with inverse $\LL'$.

\begin{lem} \label{exp} There is a natural morphism
$$
{\rm R}\Gamma_{i(P')}\Omega^t_{P \times P^*} \to \Omega_{P \times P^*}^t \underset{\rho_!\OO_{P \times P^*}}{\ltens} \rho_! \LL.
$$
\end{lem}
\dim\ \ The proof is similar to that of Lemma 3.1.2 of \cite{KS97}. We divide the proof in several steps. Set $Z=P \times P^*$ for short and let $U=\{(x,y) \in P \times P^*,\, {\rm Re}\langle x,y \rangle <-1\}$.

(i) By Lemma 3.1.2 (i) of \cite{KS97} there is a natural arrow
\begin{equation} \label{3.1.2}
\rho_!\LL' \to \rho_! \imin \rho {\rm R}\Gamma_U\ot_Z \to {\rm R}\Gamma_U\ot_Z.
\end{equation}

(ii) There is a natural isomorphism
\begin{equation} \label{3.1.1}
\Omega^t_Z \underset{\rho_!\OO_Z}{\ltens} \rho_!\LL \simeq \Omega^t_Z \underset{\rho_!\OO_Z}{\ltens} \rho_! \OO_Z(*S) \underset{\rho_!\OO_Z}{\ltens} \rho_!\LL \iso {\rm R}\Gamma_{Z \setminus S}\Omega_Z^t \underset{\rho_!\OO_Z}{\ltens} \rho_!\LL
\end{equation}
where the second isomorphism follows from the isomorphism $\tho(G,\Omega_Z) \underset{\OO_Z}{\ltens} \OO_Z(*S) \simeq \tho(G \otimes \sol(\OO_Z(*S)),\Omega_Z)$ ($G \in D^b_{\rc}(\CC_Z)$) of \cite{Bj93,KS96}.

(iii) There is a natural arrow
$$
{\rm R}\Gamma_{i(P')}\Omega^t_Z \underset{\rho_!\OO_Z}{\ltens} {\rm R}\Gamma_U\ot_Z \to {\rm R}\Gamma_{Z \setminus S}\Omega^t_Z
$$
induced by the multiplication.

Composing (i)-(iii) and using the $\rho_!\OO_Z(*S)$-module structure of ${\rm R}\Gamma_{Z \setminus S}\Omega^t_Z$ we obtain the required morphism.\\
\qed

\begin{lem} \label{proj} Let $\NN \in D^b_{q-good}(\D_{P^*})$. Set ${\cal L} \circ {\cal N} = \underline{p_1}{}_*({\cal L} \underset{\OO_{P \times P^*}}{\ltens} \imin {\underline{p_2}} \NN)$. There is a morphism
\begin{equation}\label{morproj}
\Omega^t_{P} \underset{\rho_!\D_P}{\ltens} \rho_!({\cal L} \circ \NN) \gets (\Omega^{t}_{P^*} \underset{\rho_!\D_{P^*}}{\ltens}\rho_!\NN)^\vee[-n].
\end{equation}
\end{lem}
\dim\ \ The required morphism is the composition of the following morphisms
\begin{eqnarray*}
\Omega^t_{P} \underset{\rho_!\D_P}{\ltens} \rho_!({\cal L} \circ \NN) & \simeq & Rp_{1*}(\Omega^t_{P \times P^*} \underset{\rho_!\D_{P \times P^*}}{\ltens} \rho_!({\cal L} \underset{\OO_{P \times P^*}}{\ltens} \imin {\underline{p_2}} \NN)) \\
& \simeq & Rp_{1*}((\Omega^t_{P \times P^*} \underset{\rho_!\OO_{P \times P^*}}{\ltens} \rho_!{\cal L}) \underset{\rho_!\D_{P \times P^*}}{\ltens} \rho_!\imin {\underline{p_2}} \NN) \\
& \gets & 
Rp_{1*}{\rm R}\Gamma_{i(P')}(\Omega^t_{P \times P^*} \underset{\rho_!\D_{P \times P^*}}{\ltens} \rho_!\imin {\underline{p_2}}\NN) \\
& \simeq & Rp_{1*}{\rm R}\Gamma_{i(P')}p_2^!(\Omega^t_{P^*} \underset{\rho_!\D_{P^*}}{\ltens} \rho_!\NN)[-n].
\end{eqnarray*}
The first isomorphism follows from \eqref{dirimotDmod}, the third arrow from Lemma \ref{exp} and the last isomorphism from \eqref{iminvotDmod}.\\
\qed


\begin{lem} \label{exp2} There is a natural morphism
$$
\rh_{\rho_!\OO_{P \times P^*}}(\LL,\OW_{P \times P^*|E \times P^*}) \to (\OW_{P \times P^*|P \times E^*})_{i(P)}.
$$
\end{lem}
\dim\ \ The proof is similar to that of Lemma 3.1.2 of \cite{KS97}. We divide the proof in several steps. Set $Z=P \times P^*$ for short and let $U=\{(x,y) \in P \times P^*,\, {\rm Re}\langle x,y \rangle <-1\}$.

(i) There is a natural isomorphism
\begin{eqnarray*} \label{3.1.12}
\rh_{\rho_!\OO_{Z}}(\rho_!\LL,\OW_{Z}) & \simeq & \rh_{\rho_!\OO_{Z}}( \rho_! \LL ,\rh_{\rho_!\OO_Z}( \rho_!\OO_Z(*S),\OW_Z)) \\
& \osi & \rh_{\rho_!\OO_{Z}}( \rho_! \LL ,\OW_{Z|Z \setminus S})
\end{eqnarray*}
where the second isomorphism follows from the isomorphism $\rh_{\OO_Z}(\OO_Z(*S),G \wtens \OO_Z) \simeq (\sol(\OO_Z(*S)) \otimes G) \wtens \OO_Z$ ($G \in D^b_{\rc}(\CC_Z)$) of \cite{Bj93,KS96}.

(ii) There is a natural arrow
\begin{eqnarray*}
\rh_{\rho_!\OO_{Z}}( \rho_! \LL ,\OW_{Z|Z \setminus S}) & \simeq & (\rho_!\LL' \underset{\rho_!\OO_{Z}}{\ltens} \OW_{Z|Z \setminus S})_{Z \setminus S} \\
& \to & (\rho_!\LL' \underset{\rho_!\OO_{Z}}{\ltens} \OW_{Z|Z \setminus S})_{i(P')},
\end{eqnarray*}
where the isomorphism follows from the $\rho_!\OO_Z(*S)$-module structure of $\OW_{Z|Z \setminus S}$ and the fact that $\LL'|_S=0$.

(iii) We have an arrow
\begin{eqnarray*}
(\rho_!\LL' \underset{\rho_!\OO_{Z}}{\ltens} \OW_{Z|Z \setminus S})_{i(P')} & \to & (\rho_!\imin \rho \Gamma_U\ot_{Z} \underset{\rho_!\OO_{Z}}{\ltens} \OW_{Z|Z \setminus S})_{i(P')} \\
& \to & (\OW_{Z|Z \setminus S})_{i(P')},
\end{eqnarray*}
where the first arrow follows from Lemma 3.1.2 of \cite{KS97} and the second one is induced by the multiplication.

Composing (i)-(iii) and using the fact that
there are natural arrows induced by the distinguished triangle \eqref{dtFGH}
\begin{eqnarray*}
\OW_{P \times P^*|E \times P^*} & \to & \OW_{P \times P^*}, \\
\OW_{P \times P^*|E \times E^*} & \to & \OW_{P \times P^*|P \times E^*},
\end{eqnarray*}
we obtain the required morphism.\\
\qed

\begin{lem} \label{proj} Let $\NN \in D^b_{q-good}(\D_{P^*})$. Set ${\cal L} \circ {\cal N} = \underline{p_{1}}{}_*({\cal L} \otimes_{\OO_{P \times P^*}} \imin {\underline{p_2}} \NN)$. There is a morphism
\begin{equation}\label{morproj2}
\rh_{\rho_!\D_{P}} (\rho_!({\cal L} \circ \NN),\OW_{P|E}) \to \rh_{\rho_!\D_{P^*}}(\rho_!\NN,\OW_{P^*|E^*}))^\land[n].
\end{equation}
\end{lem}
\dim\ \ The required morphism is the composition of the following morphisms
\begin{eqnarray*}
\lefteqn{\rh_{\rho_!\D_{P}} (\rho_!({\cal L} \circ \NN),\OW_{P|E})[-n]} \\
& \simeq & Rp_{1*}\rh_{\rho_!\D_{P \times P^*}}(\rho_!(\LL \underset{\OO_{P \times P^*}}{\ltens} \underline{p_2}^{-1}\NN), \OW_{P \times P^*|E \times P^*}) \\
& \simeq & Rp_{1*}\rh_{\rho_!\D_{P \times P^*}}(\rho_!\underline{p_2}^{-1}\NN,\rh_{\rho_!\OO_{P \times P^*}}(\rho_!\LL,\OW_{P \times P^*|E \times P^*})) \\
& \to & Rp_{1*}(\rh_{\rho_!\D_{P \times P^*}}(\rho_!\underline{p_2}^{-1}\NN,\OW_{P \times P^*|P \times E^*}))_{i(P')} \\
& \simeq & Rp_{1*}(\imin {p_2}\rh_{\rho_!\D_{P^*}}(\rho_!\NN,\OW_{P^*|E^*}))_{i(P')},
\end{eqnarray*}
where the first isomorphism follows from \eqref{dirimowDmodF}, the fourth arrow from Lemma \ref{exp2} and the last isomorphism from \eqref{iminvowDmodF}.\\
\qed

\begin{teo}\label{laplace} The Laplace transform induces
quasi-isomorphisms of $D(E^*_{sa})$-modules
\begin{eqnarray*}
\OO_{E_{\RP}}^{t\land}[n] & \simeq & \ot_{E_{\RP}^*}, \\
\OO_{E_{\RP}}^{\mathrm{w}\land}[n] & \simeq & \OW_{E_{\RP}^*}.
\end{eqnarray*}
\end{teo}
\dim\ \ Let $(\cdot)_{an}$ be the canonical functor sending algebraic $\D$-modules to analytic $\D$-modules on $P$. With the notations of Lemma \ref{compact}, set
$\NN=(\underline{i_2}{}_{*} \D_{E^*})_{an}$ and remark that $\LL \simeq (\underline{i}{}_*\D_{E \times E^*}e^{-\langle x,y \rangle})_{an}$. Then $\NN \circ \LL \simeq (\underline{i_1}{}_*(\D_{E^*} \circ \D_{E \times E^*}e^{-\langle x,y \rangle}))_{an} \simeq (\underline{i_1}{}_*\D_E)_{an}$. The last isomorphism is compatible with the isomorphism between Weyl algebras induced by the Fourier transform, we refer to \cite{Ma87-88} for more details. We omit the functor $(\cdot)_{an}$ to lighten notations.

(i) The Laplace transform induces a morphism of $D(E_{sa})$-modules
\begin{equation}\label{adjlaplacet}
\OO^{t\vee}_{E^*_{\RP}}[-n] \to \ot_{E_{\RP}}.
\end{equation}
First remark that \eqref{morproj} and the equivalence between left and right $\D$-modules define a morphism
$$
\rh_{\rho_!\D_{P^*}}(\rho_!\underline{i_2}{}_{*} \D_{E^*},\ot_{P^*})^\vee[-n] \to \rh_{\rho_!\D_P}(\rho_!\underline{i_1}{}_{*} \D_E,\ot_P).
$$
Then the morphism \eqref{adjlaplacet} is constructed as follows
\begin{eqnarray*}
\OO^{t\vee}_{E^*_{\RP}}[-n] & \simeq & (Rj_{2*} \rh_{\rho_!\D_{P^*}}(\rho_!\underline{i_2}{}_{*} \D_{E^*}, \OO^{t}_{P^*} ))^\vee[-n]\\
& \simeq & Rj_{1*}( \rh_{\rho_!\D_{P^*}}(\rho_!\underline{i_2}{}_{*} \D_{E^*}, \OO^{t}_{P^*} ))^\vee[-n]\\
& \to & Rj_{1*}\rh_{\rho_!\D_P}(\rho_!\underline{i_1}{}_{*} \D_E,\ot_P) \\
& \simeq & \ot_{E_{\RP}}.
\end{eqnarray*}
The second isomorphism follows from Lemma \ref{compact} and the third morphism follows from Lemma \ref{proj}.
By adjunction we obtain a morphism $\ot_{E^*_{\RP}} \to \OO^{t\land}_{E_\RP}[n]$. In a similar manner one can construct the morphism $\OW_{E^*_{\RP}} \to \OO^{{\rm w}\land}_{E_{\RP}}$.

(ii) It is enough to check the isomorphism
${\rm R}\Gamma(U;\OO_{E_{\RP}}^{t\land}[n]) \simeq
{\rm R}\Gamma(U;\ot_{E^*_{\RP}})$ on a basis for the topology of
$E_{sa}$. Hence we may assume that $U$ is $\RP$-connected and then
that
 $U$ is an open
su\-ba\-na\-ly\-tic cone of $E$. We have the isomorphisms
\begin{eqnarray*}
{\rm R}\Gamma(U;\OO_{E_{\RP}}^{t\land}[n])  & \simeq &
\Rh(\CC_U,\OO_{E_{\RP}}^{t\land}[n])\\
 & \simeq & \Rh(\CC_U^\vee[-n],\ot_{E_{\RP}})\\
 & \simeq & \Rh(\CC_U^{\land a}[n],\ot_{E_{\RP}})\\
 & \simeq & \Rh(\CC_U,\ot_{E^*_{\RP}})\\
 & \simeq & {\rm R}\Gamma(U;\ot_{E^*_{\RP}})
\end{eqnarray*}
where $a$ denotes the antipodal map. The third isomorphism follows
from Lemma 3.7.10 of \cite{KS90} and the fourth one is given by
the Laplace isomorphism of Theorem 5.2.1 of \cite{KS97}.

(iii) Similarly, let $U$ be an open su\-ba\-na\-ly\-tic l.c.t. cone of
$E$. Then $\CC_U \simeq D'\CC_{\overline{U}}$. We have the isomorphisms
\begin{eqnarray*}
{\rm R}\Gamma(U;\OO_{E_{\RP}}^{\text{w}\land}[n])  & \simeq &
\Rh(D'\CC_{\overline{U}},\OO_{E_{\RP}}^{\text{w}\land})[n]\\
 & \simeq & \Rh((D'\CC_{\overline{U}})^\vee[-n],\OW_{E_{\RP}})\\
 & \simeq & \Rh(D'(\CC_{\overline{U}}^\land[n]),\OW_{E_{\RP}})\\
 & \simeq & \CC_{\overline{U}}^\land[n] \overset{\rm W}{\otimes}\OO_E \\
 & \simeq & \CC_{\overline{U}} \overset{\rm W}{\otimes}\OO_{E^*} \\
 & \simeq & \Rh(D'\CC_{\overline{U}},\OW_{E^*_{\RP}})\\
 & \simeq & {\rm R}\Gamma(U;\OW_{E^*_{\RP}}).
\end{eqnarray*}
The  third isomorphism
follows from Proposition 3.7.12 of \cite{KS90} and the fifth one
is given by
the Laplace isomorphism of Theorem 5.2.1 of \cite{KS97}.\\

Moreover these two isomorphisms are linear over the Weyl algebra
$D(E_{sa})$.\\
\qed

\begin{oss} In \cite{KS97} the isomorphism $\ot_{E^*} \simeq \ot{}^\land_E[n]$ is established in the category $\mod_{\RP}(\CC_{E^*})$ and it is $D(E^*)$-linear (it was also later proven in \cite{Be01} with other techniques).  Applying the functor $\imin \rho$ to the isomorphism of Theorem \ref{laplace} we can recover this result.
\end{oss}




\addcontentsline{toc}{section}{

\textsc{Università di Padova, Dipartimento di Matematica Pura ed
Applicata, via Trieste 63, 35121 Padova, Italy}

\textit{E-mail address: }\verb"lprelli@math.unipd.it"


\textbf{References}}



\begin{thebibliography}{99}


 \vspace{0.5cm}


\bibitem{Be01} O. Berni; Laplace transform of temperate holomorphic functions; Compositio Mathematica $\bf{129}$ pp. 183-201 (2001).


\bibitem{BM88} E. Bierstone, D. Milmann; Semianalytic and subanalytic sets; Publ. I.H.E.S. $\bf{67}$ pp. 5-42 (1988).


\bibitem{Bj93} J. E. Bj\"ork; Analytic ${\cal D}$-modules and applications; Mathematics and its Applications, $\bf{247}$, Kluwer Academic Publishers Group, Dordrecht (1993).


\bibitem{Co00} M. Coste; An introduction to o-minimal geometry; Dip. Mat. Univ. Pisa, Dottorato di Ricerca in Matematica, Istituti Editoriali e Poligrafici Internazionali, Pisa (2000).


\bibitem{EP} M. Edmundo, L. Prelli; Sheaves on $\mathcal{T}$-topologies, arXiv:1002.0690.










\bibitem{Ka84} M. Kashiwara; The Riemann-Hilbert problem for holonomic systems; Publ. RIMS, Kyoto Univ. $\bf{20}$, pp. 319-365 (1984).


\bibitem{Ka03} M. Kashiwara; $\D$-modules and microlocal calculus; Translations of Math. Monog. $\bf{217}$, Iwanami Series in Modern Math., American Math. Soc., Providence (2003).


\bibitem{KS90} M. Kashiwara, P. Schapira; Sheaves on manifolds; Grundlehren der Math. $\bf{292}$ Springer-Verlag, Berlin (1990).


\bibitem{KS96} M. Kashiwara, P. Schapira; Moderate and formal cohomology associated with constructible sheaves; M\'emoires Soc. Math. France $\bf{64}$ (1996).


\bibitem{KS97} M. Kashiwara, P. Schapira; Integral transforms with exponential kernels and Laplace transform; Journal of the AMS $\bf{4}$ vol. $\bf{10}$ (1997).


\bibitem{KS01} M. Kashiwara, P. Schapira; Ind-sheaves; Ast{\'e}risque $\bf{271}$ (2001).




\bibitem{Lo93} S. \L ojaciewicz; Sur la géométrie semi- et sous-analytique; Ann. Inst. Fourier $\bf{43}$ pp. 1575-1595 (1993).


\bibitem{Lo59} S. \L ojaciewicz; Sur le probl\`eme de la division; Studia Mathematica $\bf{8}$ pp. 87-136 (1959).


\bibitem{Ma67} B. Malgrange; Ideals of differentiable functions; Tata Institute, Oxford University Press (1967).


\bibitem{Ma87-88} B. Malgrange; Transformation de Fourier géométrique; Séminaire Bourbaki, $\bf{30}$, Exp. N. $\bf{692}$ (1987-1988).


\bibitem{Pr1} L. Prelli; Sheaves on subanalytic sites; Rend. Sem. Mat. Univ. Padova Vol. $\bf{120}$, pp. 167-216 (2008).


\bibitem{Pr2} L. Prelli; Microlocalization of subanalytic sheaves; C. R. Acad. Sci. Paris Math. $\bf{345}$, pp. 127-132 (2007), arXiv:math.AG/0702459v3.




\bibitem{Ta94} G. Tamme; Introduction to \'etale cohomology; Universitext Springer-Verlag, Berlin (1994).


\bibitem{VD98} L. Van Der Dries; Tame topology and o-minimal structures; London Math. Society Lecture Notes Series $\bf{248}$; Cambridge University Press, Cambridge (1998).


\bibitem{Wi05} A. J. Wilkie; Covering definable open sets by open cells; in M. Edmundo, D. Richardson, and A. Wilkie, editors, O-minimal Structures, Proceedings of the RAAG Summer School Lisbon 2003, Lecture Notes in Real Algebraic and Analytic Geometry, Cuvillier Verlag (2005).




\end{thebibliography}
\end{document}